\magnification=1200
\overfullrule=0mm
\baselineskip=15pt
\centerline{\bf Linear equations over multiplicative groups, recurrences, and mixing I}
\medskip
\centerline {H. Derksen\footnote*{partially supported by the NSF, grant DMS 0901298.} and D. Masser}
\bigskip
\noindent
\bigskip
\noindent
\bigskip
\noindent
{\bf Abstract.} {\it Let $K$ be a field of positive characteristic. When $V$ is a linear variety in $K^n$ and $G$ is a finitely generated subgroup of $K^*$, we show how to compute the set $V \cap G^n$ effectively using heights. We calculate all the estimates explicitly. A special case provides the effective solution of the $S$-unit equation in $n$ variables.}
\bigskip
\noindent
{\tt 2000 MSC codes.} 11D04, 11D72, 11G35, 11G50, 14G25.
\bigskip
\bigskip
\bigskip
\noindent
{\bf 1. Introduction.} In 2004 the second author published a paper [Mass] about linear equations over multiplicative groups in positive characteristic. This was specifically aimed at an application to a problem about mixing for dynamical systems of algebraic origin, and, as a result about linear equations, it lacked some of the simplicity of the classical results in zero characteristic. A new feature was the appearance of $n-1$ independently operating Frobenius maps; here $n$ is the number of variables.
\medskip
In 2007 the first author published a paper [D] about recurrences in positive characteristic. He proved an analogue of the famous Skolem-Lech-Mahler Theorem in zero characteristic. A new feature was the appearance of integer sequences involving combinations of $d-2$ powers of the characteristic; here $d$ is the order of the recurrence.
\medskip
It turns out that these two new features are identical. In positive characteristic the vanishing of a recurrence with $d$ terms can be regarded as an linear equation in $d-1$ variables to be solved in a multiplicative group (so in particular $n-1=d-2$). This observation will be developed in three directions.
\medskip
In the present paper we give an improved version of the result of [Mass] in a form more closely related to that in zero characteristic. In fact we shall prove some quantitative versions in which all the estimates are effective and furthermore we shall make them completely explicit. This is in sharp contrast to the situation in zero characteristic, where even in very simple circumstances there are no effective upper bounds for the solutions.
\medskip
In a second paper we shall apply these results to recover the main theorem of [D], which we even generalize to sums of recurrences. In zero characteristic rather little is known about such sums, and indeed there is a conjecture of Cerlienco, Mignotte and Piras [CMP] to the effect that such problems are undecidable. In positive characteristic we will establish not only the decidability but also give completely effective algorithms to solve the problem. 
\medskip
In a third paper we apply our linear equations results to give an algorithm to determine the smallest order of non-mixing of any basic action associated with a given prime ideal in a Laurent polynomial ring. From [Mass] we know that the non-mixing comes from the so-called non-mxing sets, and our work even provides a way of finding these. Again the algorithms are completely effective.
\medskip
We begin by recalling the classical result for a linear equation in zero characteristic, for convenience in homogeneous form. For a field $K$ we write $K^*$ for the multiplicative group of all non-zero elements of $K$. For any subgroup $G$ of $K^*$ and a positive integer $n$ it makes sense to write ${\bf P}_n(G)$ for the set of points in projective space defined over $G$.
\bigskip
\noindent
{\bf Theorem A} (Evertse [E], van der Poorten-Schlickewei [PS]).  {\it Let $K$ be a field of zero characteristic, and for $n \geq 2$ let $a_0,\ldots,a_n$ be non-zero elements of $K$. Then for any finitely generated subgroup $G$ of $K^*$ the equation 
$$a_0X_0+a_1X_1+\cdots+a_nX_n~=~0 \eqno(1.1)$$
has only finitely many solutions $(X_0,X_1,\ldots,X_n)$ in ${\bf P}_n(G)$ which satisfy
$$\sum_{i \in I}a_iX_i \neq 0 \eqno(1.2)$$
for every non-empty proper subset $I$ of $\{0,1,\ldots,n\}$.}
\bigskip
We should point out that this remains true even when $G$ is not finitely generated but has finite ${\bf Q}$-dimension. See also a recent paper [EZ] of Evertse and Zannier for an interesting function field version.
\medskip
Theorem A is false in positive characteristic $p$; for example in inhomogeneous form for $n=2$ the equation $$x+y=1 \eqno(1.3)$$ has a solution $x=t,~y=1-t$ over the group $G$ in $K={\bf F}_p(t)$ generated by $t,1-t$; and so thanks to Frobenius infinitely many solutions 
$$x=t^{p^e},~y=1-t^{p^e}=(1-t)^{p^e}~~~(e=0,1,2,\ldots) \eqno(1.4)$$ 
which all satisfy (1.2).
\medskip
We can regard Theorem A as a descent step from the hyperplane $H$ defined by equation (1.1) to proper linear subvarieties defined by the vanishing of the left-hand sides in (1.2). We can iterate this descent by introducing special varieties $T$ defined solely by binary equations of the shape $X_i=aX_j~(i \neq j,a \neq 0)$. For example $T$ could be a single point or, when there are no equations at all, the full ${\bf P}_n$. We could call such varieties linear cosets or just cosets. This word has a group-theoretical connotation, and indeed $T$ above is a translate of a group subvariety of the multiplicative group ${\bf G}_{\rm m}^n$ in ${\bf P}_n$. Conversely it is not difficult to see that every linear translate of a group subvariety of ${\bf G}_{\rm m}^n$ is a coset in our sense (see for example Lemma 9.4 p.76 of [BMZ]). But we will in this paper make no use of these remarks or indeed hardly any further reference to group varieties.
\medskip
Anyway, it is easily seen that the complete descent yields a finite collection of cosets $T$, each contained in the original $H$, such that the full solution set $H(G)=H \cap {\bf P}_n(G)$ coincides with the union of all $T(G)=T \cap {\bf P}_n(G)$. This is a little closer to the more general context of Mordell-Lang (see below). No further descent from $T(G)$ in terms of proper subvarieties is possible; by way of compensation it is very simple to describe $T(G)$ explicitly (see for example the discussion towards the end of section 12).
\medskip
In positive characteristic we can establish a descent step similar to Theorem A, but it may involve Frobenius as in (1.4). This less simple situation makes the iteration more problematic, and for this reason it is clearer to present our result as a descent now from an arbitrary linear variety $V$ to proper linear subvarieties.
\medskip
However the Frobenius does not always generate infinitely many solutions. It does above for $x+y=1$, and also for 
$$t^mx+y=1 \eqno(1.5)$$
by taking a new variable $t^mx$; this is because $t$ lies in $G$. The situation is slightly more subtle for (1.5) over the group $G_l$ generated by $t^l$ and $1-t$; the above solution of (1.3) certainly leads to solutions
$$x=t^{-m}t^{p^e},~y=(1-t)^{p^e}~~~(e=0,1,2,\ldots),\eqno(1.6)$$ 
but these will not be over $G_l$ unless $p^e \equiv m$ mod $l$. This can however happen for infinitely many $e$ but not necessarily all $e$ in (1.6). This time $t$ may not lie in $G_l$ but some positive power does. Finally the equation $(1+t)x+y=1$ has a solution $x=1-t,y=t^2$ over $G$, but the use of Frobenius will bring in an extra $1+t$, no positive power of which is in $G$ (provided $p \neq 2$).
\medskip  
These considerations lead naturally to the radical $\sqrt{G}=\root K \of G$ for general $G$ in general $K^*$. For us this remains in $K$; thus it is the set of $\gamma$ in $K$ for which there exists a positive integer $s$ such that $\gamma^s$ lies in $G$. Usually $K$ will be finitely generated over its prime field, and then it is well-known that the finite generation of $G$ is equivalent to that of $\sqrt{G}$. We also see the need for some concept of isotriviality, already present in diophantine geometry at least since N\'eron's 1952 proof of the relative Mordell-Weil Theorem and Manin's 1963 proof of the relative Mordell Conjecture. In our linear context the appropriate refinement is $G$-isotriviality, introduced by Voloch [V] for $n=2$. 
\medskip
Namely, let $K$ be a field of positive characteristic $p$, and for $n \geq 2$ let $V$ be a linear variety in ${\bf P}_n$ defined over $K$. We say that $V$ is $G$-isotrivial if there is an automorphism $\psi$ of ${\bf P}_n(K)$, defined by
$$\psi(X_0,\ldots,X_n)~=~(g_0X_0,\ldots,g_nX_n) \eqno(1.7)$$
with $g_0,\ldots,g_n$ in $G$, such that $\psi(V)$ is defined over the algebraic closure $\overline {{\bf F}_p}$. Such a $\psi$ could be called a $G$-automorphism. Let us write ${\bf F}_K$ for $\overline {{\bf F}_p} \cap K$; then of course $\psi(V)$ is defined over ${\bf F}_K$. So $\psi(V)$ is defined over some ${\bf F}_q$; and now a point $w$ on $V$ defined over $G$ gives $\psi(w)$ on $\psi(V)$ which by Frobenius leads to points $\psi(w)^{q^e}~(e=0,1,2,\ldots)$ on $\psi(V)$ and so 
$$\psi^{-1}(\psi(w)^{q^e})~~~~(e=0,1,2,\ldots)\eqno(1.8)$$ 
on $V$, all still defined over $G$.
\medskip
Of course points over $G$ are nothing other than zero-dimensional $G$-isotrivial varieties.
\medskip
Here is a preliminary version of our main descent step on linear equations. For $V$ as above write $V(G)=V \cap {\bf P}_n(G)$ for the set of points of $V$ defined over $G$. But it is clearer first to consider points over the radical $\sqrt{G}$. 
\bigskip
\noindent
{\bf Descent Step over $\sqrt{G}$}. {\it Let $K$ be a field of positive characteristic, and suppose that the positive-dimensional linear variety $V_0$ defined over $K$ is not a coset. Suppose also that $\sqrt{G}$ in $K$ is finitely generated. Then there is an effectively computable finite collection ${\cal W}$ of proper $\sqrt{G}$-isotrivial linear subvarieties $W$ of $V_0$, also defined over $K$, with the following property.

\noindent
(a) If $V_0$ is not $\sqrt{G}$-isotrivial, then
$$V_0(\sqrt{G})~=~\bigcup_{W \in {\cal W}}W(\sqrt{G}).$$

\noindent
(b) If $V_0$ is $\sqrt{G}$-isotrivial and $\psi(V_0)$ is defined over ${\bf F}_q$, then
$$V_0(\sqrt{G})~=~\psi^{-1}\left(\bigcup_{W \in {\cal W}}~\bigcup_{e=0}^\infty(\psi(W)(\sqrt{G}))^{q^e}\right).$$}
\bigskip
Thus (a) says that the points of $V_0(\sqrt{G})$ are not Zariski-dense in $V_0$; and (b) says that the points on $V_0(\sqrt{G})$ like (1.8), which can be dense, at least arise from a set of $w$ which is not dense.
\medskip
Part (a) was essentially proved for $n=2$ as Theorem 1 by Voloch [V] (p.196), and his Theorem 2 (p.198) even covers the more general case of finite ${\bf Q}$-dimension; here one gets the finiteness of the solution set. A forerunner of part (b) for $n=2$ can be seen in Mason [Maso] (pp. 107,108). The main result of [Mass] is restricted to a single equation (1.1) and is expressed in terms of a concept of ``broad" set; as we do not need this result here (or even the concept) we refrain from quoting it. However these authors do not discuss the effectivity in our sense (see the discussion below).
\medskip
A simple example of (b) in inhomogeneous form is (1.3); this represents a line $L$, clearly isotrivial and even trivial in that we can take $\psi$ as the identity automorphism. When $G$ is generated by $t$ and $1-t$ in $K={\bf F}_p(t)$, then $\sqrt{G}$ is obtained by adding the elements of ${\bf F}_p^*$ as generators. Leitner [Le] has found that for $p \geq 3$ there are $p+4$ points $W$, six of which are like $w=(t,1-t)$ in (1.4) and the remaining $p-2$ are the $w=(x,1-x)$ for $x=2,3,\ldots,p-1$.
\medskip
So much for $V_0(\sqrt{G})$. In the analogous characterization of $V_0(G)$ there is no longer a clear separation of cases. In fact it can happen in case (b) above that the actions of Frobenius through $q^e$ can get truncated, so that each $e$ remains bounded; but then it is easy to reduce this to case (a). A simple example is (1.5) for $m=1$ in the group $G=G_l$ above for $l=p$, when the solutions (1.6) are over $G$ only when $e=0$. Here is a general statement.
\bigskip
\noindent
{\bf Descent Step over $G$}. {\it Let $K$ be a field of positive characteristic, and suppose that the positive-dimensional linear variety $V_0$ defined over $K$ is not a coset. Suppose also that $\sqrt{G}$ in $K$ is finitely generated. Then there is an effectively computable finite collection ${\cal W}$ of proper $\sqrt{G}$-isotrivial linear subvarieties $W$ of $V_0$, also defined over $K$, such that either
$$V_0(G)~=~\bigcup_{W \in {\cal W}}W(G)$$
or
$$V_0(G)~=~\psi^{-1}\left(\bigcup_{W \in {\cal W}}~\bigcup_{e=0}^\infty(\psi(W)(G))^{q^e}\right)$$
for some $q$ and some $\sqrt{G}$-automorphism $\psi$ with $\psi(V_0)$ defined over ${\bf F}_q$.}
\bigskip
\medskip
It may be instructive here to consider the inhomogeneous example
$$x+y-z~=~1\eqno(1.9)$$
still over the group $G$ in $K={\bf F}_p(t)$ generated by $t,1-t$. Now (1.9) represents a plane $P$, also isotrivial and even trivial. Leitner [Le] has found that for $p \geq 5$ there are 22 lines $W$ and 8 points $W$. For example the line defined by 
$$tx+y=1,~~~z=(1-t)x \eqno(1.10)$$
is one of these. So is the coset line defined by $x=z,~y=1$. And so is the point
$$x=t,~~y={1-t \over t},~~z={(1-t)^2 \over t}.$$
\medskip
We can easily iterate the descent from (1.10). This is isotrivial via the automorphism $\psi$ taking $x,y,z$ to $\tilde x=tx,\tilde y=y,\tilde z={t \over 1-t}z$, when the equations become $\tilde x+\tilde y=1,~\tilde z=\tilde x$. Now (1.4) (with $e$ replaced by $f$) on (1.3) lead to the points $w=(x,y,z)$ of $W(G)$ with
$$x=t^{p^f-1},~~~y=(1-t)^{p^f},~~~z=t^{p^f-1}(1-t)~~~(f=0,1,2,\ldots).$$
Then from (1.8) (with $q=p$ and the identity automorphism) we get the points
$$x=t^{(q-1)r},~~~y=(1-t)^{qr},~~~z=t^{(q-1)r}(1-t)^r\eqno(1.11)$$
of $P(G)$; here $q=p^f$ and $r=p^e$ now indicate independently varying powers of $p$. This is precisely the example in [Mass] (p.202).
\medskip
With the help of a suitable notation we can after all do the complete descent, also for linear varieties that are cosets; then the latter arise solely as obstacles. Denote by $\varphi=\varphi_q$ the Frobenius with $\varphi(x)=x^q$. Let $\psi_1,\ldots,\psi_h$ be projective automorphisms. Then we imitate commutator brackets by defining the operator 
$$[\psi_1,\ldots,\psi_h]~~=~~[\psi_1,\ldots,\psi_h]_q~~=~~\bigcup_{e_1=0}^\infty \cdots \bigcup_{e_h=0}^\infty (\psi_1^{-1}\varphi^{e_1} \psi_1) \cdots (\psi_h^{-1}\varphi^{e_h} \psi_h),\eqno(1.12)$$
with of course the identity interpretation if $h=0$. This formally resembles Definition 7.7 of [D] (p.208).
\bigskip
\noindent
{\bf Theorem 1.} {\it Let $K$ be a field of positive characteristic $p$, let $V$ be an arbitrary linear variety defined over $K$, and suppose that $\sqrt{G}$ in $K$ is finitely generated. Then there is a power $q$ of $p$ such that $V(G)$ is an effectively computable finite union of sets $[\psi_1,\ldots,\psi_h]_qT(G)$ with $\sqrt{G}$-automorphisms $\psi_1,\ldots,\psi_h ~(0 \leq h \leq n-1)$, and cosets $T$ contained in V.}
\bigskip
Here we see quite clearly the $n-1$ Frobenius operators mentioned in the first paragraph of section 1. In general they act independently because they are separated by automorphisms. The example
$$x_1+x_2-x_3-\cdots-x_n~~=~~1$$ 
generalizes (1.3) and (1.9), and it can be used to show that the upper bound $n-1$ in Theorem 1 cannot always be improved. This we carry out in section 13 on limitation results. The same can also be seen indirectly through the applications to recurrences, where we will see that the analogous upper bound $d-2$ cannot always be improved. 
\medskip
Taking $e_1=1$ in (1.12) and all other zero, we see that $\psi_1^{q-1}$ is a $G$-automorphism. Similarly for $\psi_1^{q-1},\ldots,\psi_h^{q-1}$. However it may not always be possible to choose $\psi_1,\ldots,\psi_h$ as $G$-automorphisms. This we also prove in section 13. 
\medskip
We can also symmetrize the sets in Theorem 1. We explain this with the points (1.11) on $P$ defined by (1.9). They can be written as
$$x=t^{s-r},~~~y=(1-t)^s,~~~z=t^{s-r}(1-t)^r\eqno(1.13)$$
with $s=qr$. Here there is asymmetry because apparently $r$ divides $s$. However (1.13) has a meaning for any independent positive powers $r,s$ of $p$; and it is easily checked that the resulting points remain on $P$.
\medskip
To formulate this in general we introduce another bracket notation more related to the group law. For points $\pi_0,\pi_1,\ldots,\pi_h$ we define the set
$$(\pi_0,\pi_1,\ldots,\pi_h)~~=~~(\pi_0,\pi_1,\ldots,\pi_h)_q~~=~~\pi_0\bigcup_{l_1=0}^\infty \cdots \bigcup_{l_h=0}^\infty (\varphi^{l_1} \pi_1) \cdots (\varphi^{l_h} \pi_h), \eqno(1.14)$$
with of course the interpretation $\pi_0$ itself if $h=0$. We introduce more special varieties $S$ defined solely by binary equations of the shape $X_i=X_j$. For example $S$ could be the single point with all coordinates equal or the full ${\bf P}_n$. We could call such varieties linear subgroups or just subgroups. As above it is not difficult to see that they are precisely the linear group subvarieties of ${\bf G}_{\rm m}^n$, but again we don't need to know this.
\bigskip
\noindent
{\bf Theorem 2.} {\it Let $K$ be a field of positive characteristic $p$, let $V$ be an arbitrary linear variety defined over $K$, and suppose that $\sqrt{G}$ in $K$ is finitely generated. Then there is a power $q$ of $p$ such that $V(G)$ is an effectively computable finite union of sets $(\pi_0,\pi_1,\ldots,\pi_h)_qS(G)$ with points $\pi_0,\pi_1,\ldots,\pi_h ~(0 \leq h \leq n-1)$ defined over $\sqrt{G}$ and subgroups $S$.}
\bigskip
As in Theorem 1, the upper bound $n-1$ in Theorem 2 cannot always be improved. We shall verify this in section 13. Also one can easily see that $\pi_0^{q-1},\pi_1^{q-1},\ldots,\pi_h^{q-1}$ (as well as the product $\pi_0\pi_1\cdots\pi_h$) are defined over $G$. However this may not always be true of $\pi_0,\pi_1,\ldots,\pi_h$, as we shall also prove in section 13.
\medskip
When $V$ is a hyperplane defined by (1.1) we can even descend to points, provided we restrict to (1.2) in the style of Theorem A.
\bigskip
\noindent
{\bf Theorem 3.} {\it Let $K$ be a field of positive characteristic $p$, let $H$ be defined by 
$$a_0X_0+a_1X_1+\cdots+a_nX_n~=~0$$ 
for non-zero $a_0,a_1,\ldots,a_n$ in $K$, and write $H^*(G)$ for the set of points in ${\bf P}_n(G)$ satisfying 
$$\sum_{i \in I}a_iX_i \neq 0$$ 
for every non-empty proper subset $I$ of $\{0,1,\ldots,n\}$. Suppose that $\sqrt{G}$ in $K$ is finitely generated. Then there is a power $q$ of $p$ such that $H^*(G)$ is contained both in (1) an effectively computable finite union of sets $[\psi_1,\ldots,\psi_h]_q\{\tau\}$ in $H(G)$ with $\sqrt{G}$-automorphisms $\psi_1,\ldots,\psi_h ~(0 \leq h \leq n-1)$ and points $\tau$, and in (2) an effectively computable finite union of sets $(\pi_0,\pi_1,\ldots,\pi_h)_q$ in $H(G)$ with points $\pi_0,\pi_1,\ldots,\pi_h ~(0 \leq h \leq n-1)$.}
\bigskip
We do not prove it here, but in this situation $H^*(G)$ is precisely a finite union of $[\psi_1,\ldots,\psi_h]_q\{\tau\}$. However there seems to be a strange asymmetry between the asymmetric part (1) and the symmetric part (2). Namely it seems improbable that  $H^*(G)$ is precisely a finite union of $(\pi_0,\pi_1,\ldots,\pi_h)_q$. For example, the point (1.13) on $H$ defined by (1.9) is in $H^*(G)$ except for $r=s$, which disturbs the independence of $r$ and $s$.
\medskip
Apart from the work [V] already mentioned, there are other results of this kind, now in the more general context of Mordell-Lang for arbitrary varieties $V$ inside arbitrary semiabelian varieties $\bf S$. Typically here one intersects $V$ with a finitely generated subgroup $\Gamma$ of $\bf S$; however in the present paper with ${\bf S}={\bf G}_{\rm m}^n$ we have for simplicity restricted $\Gamma$ to a Cartesian product $G^n$.
\medskip
Thus the main result Theorem A (p.104 - see also p.109) of Abramovich and Voloch [AV] almost implies part (a) of our Descent Step over $\sqrt{G}$, except that they assume that $V$ is not $K^*$-isotrivial and they have no information about $W$ which would ensure linearity in our situation. The main result Theorem 1.1 (p.667) of Hrushovsky's well-known paper [Hr] gives a similar implication. The restriction to our (a) corresponds to their restriction to the non-isotrivial case. Again these authors do not discuss the effectivity in our sense.
\medskip
After earlier work by Scanlon, the isotrivial case was treated by Moosa and Scanlon.  Their Theorem B (p.477) of [MS2] implies that our $V(G)$ is what they call an $F$-set (see also [MS1]). Indeed in our situation and notation an $F$-set is nothing other than a finite union of $(\pi_0,\pi_1,\ldots,\pi_h)_qA(G)$ with $\pi_0\pi_1\cdots \pi_h$ and $\pi_0^{q-1},\pi_1^{q-1},\ldots,\pi_h^{q-1}$ defined over $G$ and an algebraic subgroup $A$. However they do not prove the bound $h \leq n-1$ and they do not give an estimate for $A$ which would imply that it is linear because our $V$ is. Their ideas were developed by Ghioca [Gh], who in addition extended the results to Drinfeld modules. See also the work [GM] of Ghioca and Moosa on division groups. Again there is no mention of effectivity.
\medskip
Now let us discuss this effectivity, a key aspect of the present paper.
\medskip
It is well-known that Theorem A (in zero characteristic) is semieffective in the sense that effective and even explicit upper bounds for the number of solutions of (1.1) subject to (1.2) can be found. However it is not fully effective in the sense that no upper bounds are known for the size of the solutions, even in very simple cases like $K={\bf Q}$ and $G$ generated by 3,5,7; and it is even unknown how to find all the finitely many non-negative integers $a,b,c$ satisfying an equation like
$$3^a+5^b-7^c~=~1.$$
\medskip
Out of the works in positive characteristic quoted above, only two discuss effectivity, and then only semieffectivity in the sense above. Voloch [V] in the theorems mentioned above gives explicit upper bounds for the cardinality of $V(G)$ for $n=2$ in case (a) of Theorem 1; these are uniform in the sense that they are independent of $V$ and further they depend on $G$ only with regard to its rank. A similarly uniform bound is given as Theorem 6.1 (p.687) by Hrushovsky [Hr] for $V$ in an abelian variety; however as it stands it is not completely explicit due to the use of non-standard analysis. These bounds are in line with the well-known estimates in zero characteristic - see for example Theorem 1.1 of [ESS] (p.808).
\medskip
By contrast our results above are fully effective. This should be no surprise; for example it is rather easy by differentiating to find all non-negative integers $a,b,c$ with
$$(3+t)^a+(5+t)^b-(7+t)^c~=~1$$
in any fixed $K={\bf F}_p(t)$. We shall work out explicit bounds, at first for the Descent Step over $\sqrt{G}$, where the exponents appearing can be reasonably small; and then for the Descent Step over $G$ and Theorems 1,2 and 3. It would then be a straightforward matter to deduce bounds for the various cardinalities involved; but more work may be needed to make these uniform in the sense above. 
\medskip
In fact the size bounds cannot be uniform in this sense. For example from the non-isotrivial equation $x+ay=1$ with $a={1-t^m \over (1-t)^m}~(m \neq p^e)$ over the group generated by $t$ and $1-t$ in ${\bf F}_p(t)$, with solution $x=t^m,~y=(1-t)^m$, we can easily show that the size of solutions for fixed $G$ must depend on $V$. Similarly the isotrivial equation $x+y=1$ over the group generated by $t^m$ and $(1-t)^m$ in ${\bf F}_p(t)$, with the same solution, demonstrates that the size of solutions for fixed $V$ must depend on more than just the rank of $G$. 
\medskip
Because all our varieties are linear, we can measure them in a traditional way in terms of certain heights on the Grassmannian. We will show for example in the Descent Step over $\sqrt{G}$ that 
$$h(W) \leq Ch(V_0)^{2n} \eqno(1.15)$$ 
if $W$ is no longer required to be $\sqrt{G}$-isotrivial, where $C$ depends only on $K,n$ and $G$. If we insist on $W$ being $\sqrt{G}$-isotrivial, then the exponent is not so small. The well-known Northcott Property of heights often implies that the set of $W$ in (1.15) is finite and easily effectively computable.
\medskip
Perhaps since the results in zero characteristic are not effective, there is no tradition about measuring the groups $\Gamma$, even in ${\bf S}={\bf G}_{\rm m}^n$. Because our $\Gamma=G^n$, it is here possible to use a basis-free notion of regulator $R(G)$. We will show that the bounds, at least when $G=\sqrt{G}$, are all of polynomial growth in $R(G)$. For example in (1.15) we get
$$C \leq cR(G)^{6n+2}$$
again if $W$ is no longer required to be $\sqrt{G}$-isotrivial, where $c$ now depends only on $K,n$ and the rank $r$ of $G$. In fact here
$$c~=~8n^2d(10n^3(n+r)^{3(n+r)})^{2n+1}$$
with $d$ depending only mildly on $K$; for example $d=1$ if $K$ is a field of rational functions in several independent variables over a finite field.
\medskip
However we did find it a small surprise to discover that when $G \neq \sqrt{G}$ the smallest bounds can be exponential in $R(G)$. A hint of this can be seen from the above discussion of (1.5) and $G_l$. For example the simplest solution of the equation
$$t^{42}x+y=1$$
with $x,y$ in the group generated by $t^{83}$ and $1-t$ in ${\bf F}_2(t)$ is
$$x=(t^{83})^{29130742641316365655570},~~y=(1-t)^{2417851639229258349412352}; \eqno(1.16)$$
while the regulator is only $83\sqrt{3}$. For an explanation see the end of section 11.
\medskip
In section 12 we estimate the heights (in a natural sense) of all the quantities occurring in our Theorems. The bounds are polynomial in $h(V)$ and $R(G)$ if $G=\sqrt{G}$; but otherwise they may involve an extra, possibly unavoidable, exponential dependence on $R(G)$. Here too there is a Northcott Property to ensure effectivity.
\medskip
At first sight it may seem that the methods of [Mass] and [D] are unrelated. But there are close connections, and we give some hints of this in our exposition. Here we mention just that [Mass] works with derivatives and [D] works with $p$-automata and ``free Frobenius splitting". For example over ${\bf F}_p(t)$, [Mass] (p.196) has $\delta_i=({d \over dt})^i~(i=0,\ldots,p-1)$ while [D] (p.198) splits ${\bf F}_p(t)$ into a direct sum of one-dimensional ${\bf F}_p(t^p)$-subspaces $V_i~(i=0,\ldots,p-1)$ and considers the associated projections $\lambda_i$. In the natural case $V_i=t^i{\bf F}_p(t^p)$ one checks easily that the vectors $(\delta_0,t\delta_1\ldots,t^{p-1}\delta_{p-1})$ and $(\lambda_0,\lambda_1,\ldots,\lambda_{p-1})$ are connected via an invertible matrix over ${\bf F}_p$. So in some sense differentiating is equivalent to projecting. We can also quote Hrushovsky [Hr] (p.669) {\it ``Distinguishing a basis for $K/K^p$ has the effect of fixing also a stack of Hasse derivations."} As a matter of fact we do not use Hasse derivations in this paper (see the remarks at the end of section 5).
\medskip
Here is a brief section-by-section account of what follows. 
\medskip
We begin in section 2 by explaining heights. Then in section 3 we introduce derivations, and we use all this to give preliminary effective versions of the two main technical results of [Mass] about dependence over the field of differential constants.
\medskip
In section 4 we explain regulators, and in section 5 we use these to refine the work of section 3.
\medskip
Then section 6 contains a technical result which enables us to identify isotriviality, and in section 7 we record some observations about automorphisms and heights of varieties $V$.
\medskip
We are now in a position in section 8 to make effective the main argument of [Mass] yielding the subvarieties $W$, at least for points over $\sqrt{G}$ and when $V$ is either a hyperplane or trivial. We treat general $V$ in section 9 but omitting the isotriviality of the $W$. This omission is then remedied in section 10 with a simple inductive argument, and in section 11 we show how to treat points over $G$. We can then in section 12 prove effective versions of our Descent Steps and Theorems.
\medskip
Finally in section 13, as already mentioned, we show that various aspects of our results cannot be further improved.
\medskip
We would also like to draw attention to a very recent manuscript [AB] of Adamczewski and Bell for further work in the context of $p$-automata; in particular this covers also equations (1.1) and recurrences.
\bigskip 
\bigskip
\noindent
{\bf 2. Heights.} The Theorems above for arbitrary fields can easily be reduced to the case when the field is finitely generated over its ground field ${\bf F}_p$ (see section 12 below). In general let $K$ be finitely generated over a subfield $k$ in any characteristic. We shall define heights on $K$ relative to $k$; thus we suppose that $K$ is a transcendental extension of $k$. Here we do not know any basis-free notion of height, and thus we choose a transcendence basis ${\cal B}$ of $K$ over $k$ with elements $t_1,\ldots,t_b$ regarded as independent variables over $k$. The height $\tilde h(a)=\tilde h_{\cal B}(a)$ of an element $a \neq 0$ of $k[{\cal B}]=k[t_1,\ldots,t_b]$ will be its total degree $\deg a$ regarded as a polynomial; also $\tilde h(0)=0$. The height can be extended to an element $x$ of the quotient field $k({\cal B})=k(t_1,\ldots,t_b)$ by writing $x={a_1 \over a_0}$ for coprime polynomials $a_0,a_1$ in $k[{\cal B}]$ and defining
$$\tilde h(x)~=~\tilde h_{\cal B}(x)~=~\max\{\deg a_0,\deg a_1\}. \eqno(2.1)$$
That suffices for most examples, but for mixing problems we have to extend further to all of $K$. This is a standard matter using valuations.
\medskip
There is a valuation on $k[{\cal B}]$ corresponding to total degree and defined by $|a|_\infty=\exp(\deg a) ~(a \neq 0)$; and of course $|0|_\infty=0$. This extends at once to $k({\cal B})$ by multiplicativity. And for every irreducible $p$ in $k[{\cal B}]$ there is a valuation defined on $k[{\cal B}]$ by $|a|_p=\exp(-\omega_p(a)\deg p)~(a \neq 0),$ where $\omega_p(a)$ is the exact power of $p$ dividing $a$; and again $|0|_\infty=0$. And it too extends to $k({\cal B})$ by multiplicativity. Using $v$ to run over $\infty$ and all the $p$, we have the product formula $\prod_v |x|_v=1~(x \neq 0)$ and the height formula $\tilde h(x)=\log \prod_v \max\{1,|x|_v\}$.
\medskip
Now $K$ is a finite extension of $k({\cal B})$, say of degree $d$. Thus each valuation $v$ has finitely many extensions $w$ to $K$, written $w|v$. In fact
$$|x|_w=|N(x)|_v^{1/d_w}, \eqno(2.2)$$
where the norm is from the completion $K_w$ to the completion $k({\cal B})_v$ and $d_w$ is the relative degree. We also have $\sum_{w|v}d_w=d$. Now the product formula
$$\prod_w |x|_w^{d_w}=1~~~(x \neq 0)$$
holds. Further the formula
$$\tilde h(x)~=~{1 \over d}\log \prod_w \max\{1,|x|_w^{d_w}\}$$
extends the height $\tilde h=\tilde h_{\cal B}$ to an absolute height on $K$. For all this see [La2] (pp.1-19) or [BG] (pp.1-10). 
\medskip
Actually for convenience in estimating we will use from now on the relative height 
$$h(x)=h_{\cal B}(x)=d\tilde h(x) \geq 1.$$
This can be calculated directly from the minimum polynomial in the following extension of (2.1).
\bigskip
\noindent
{\bf Lemma 2.1.} {\it Suppose $x$ in $K$ satisfies an equation $A(x)=0$ with $A(t)=a_0t^e+\cdots+a_e$ for  $a_0,\ldots,a_e$ in $k[{\cal B}]$ and $A(t)$ irreducible over $k[{\cal B}]$. Then $eh(x)=d\max\{\deg a_0,\ldots,\deg a_e\}.$}
\bigskip
\noindent
{\it Proof.} Over a splitting field $L$ we have $A(t)=a_0(t-x_1)\cdots(t-x_e)$, and we can extend, keeping the same notation, all the valuations to $L$. Then Gauss's Lemma gives
$$\max\{|a_0|_w,\ldots,|a_e|_w\}~=~|a_0|_w\max\{1,|x_1|_w\}\cdots\max\{1,|x_e|_w\}.$$
If $w$ does not divide $\infty$ then the left-hand side is 1 because $a_0,\ldots,a_e$ are coprime; and otherwise they are all $\max\{|a_0|_{\infty},\ldots,|a_e|_{\infty}\}$. Taking the product with exponents $d_w$ and then taking logarithms gives on the left-hand side $d\max\{\deg a_0,\ldots,\deg a_e\}$ and on the right-hand side $h(x_1)+\cdots+h(x_e)$. This last is just $eh(x)$ because $x_1,\ldots,x_e$ are conjugate over $k({\cal B})$.
\bigskip
An immediate consequence of Lemma 2.1 is the Northcott Property; namely that for any $H$ there are at most finitely many $x$ in $K$ with $h(x) \leq H$. 
\medskip
We will also need the standard extensions to vectors. So for $x_1,\ldots,x_l$ in $K$ we define
$$h(x_1,\ldots,x_l)~=~\log \prod_w \max\{1,|x_1|_w^{d_w},\ldots,|x_l|_w^{d_w}\}.$$
For example $h(a_0,\ldots,a_e)$ in the situation of Lemma 2.1 is just $d\max\{\deg a_0,\ldots,\deg a_e\}.$ The Northcott Property extends at once to $K^l$.
\bigskip 
\bigskip
\noindent
{\bf 3. Dependence with heights.} Given $K$ finitely generated and transcendental over $k$, there is always a separable transcendence basis ${\cal B}=(t_1,\ldots,t_b)$; this means that $K$ is separable over $k({\cal B})$. As above write $d=[K:k({\cal B})].$ On $k[{\cal B}]$ we have the standard derivations ${\partial \over \partial t_1},\ldots,{\partial \over \partial t_b}$, which extend in the obvious way to $k({\cal B})$. And by separability they extend uniquely to $K$. For all this see [La1] (pp.183-184). For an integer $i \geq 0$ we define ${\cal D}(i)$ as the set of operators
$$D~=~\left({\partial \over \partial t_1}\right)^{i_1}\cdots\left({\partial \over \partial t_b}\right)^{i_b}$$
as $i_1,\ldots,i_b$ run over all non-negative integers with $i_1+\cdots+i_b \leq i$. This is not quite the same as [Mass] (p.196), where we had $i \geq 1$ and $i_1+\cdots+i_b < i$.
\medskip
It will be convenient for later calculations to define a quantity $h(x;i)$ as follows. We order in some way the operators $D_1,\ldots,D_l$ of ${\cal D}(i)$, and we define for $x \neq 0$
$$h(x;i)~=~h_{\cal B}(x;i)~=~h\left({D_1x \over x},\ldots,{D_lx \over x}\right)$$
of course independent of the ordering.
\medskip
The next result is an explicit version of Lemma 3 of [Mass] (p.195) however without reference to any group $G$. We write $C$ for the field of differential constants in $K$. For zero characteristic this is $k$, but for positive characteristic $p$ it is the set of $p$th powers of elements of $K$.
\bigskip
\noindent
{\bf Lemma 3.1.} {\it For $m \geq 2$ suppose $c_1,\ldots,c_m$ are in $C$ and $x_1,\ldots x_m$ are in $K^*$ with 
$$c_1x_1+ \cdots+c_mx_m~=~1. \eqno(3.1)$$
Then either

(a) $h(c_1x_1,\ldots,c_mx_m) ~\leq~  (m+1)\left(h(x_1;m-1)+\cdots+h(x_m;m-1)\right)$

\noindent
or 

(b) $x_1,\ldots,x_m$ are linearly dependent over $C$.}
\bigskip
\noindent
{\it Proof.} If (b) does not hold, then the theory of the generalized Wronskian (see for example [La2] p.174) shows that we may find operators $D_i$ in ${\cal D}(i) ~(i=0,\ldots,m-1)$ such that the matrix with entries $D_ix_j ~(i=0,\ldots,m-1;~j=1,\ldots,m)$ is non-singular. Applying them to (3.1) we get 
$$\sum_{j=1}^m{D_ix_j \over x_j}(c_jx_j)~~=~~D_i(1) ~~~(i=0,\ldots,m-1).$$
These can be solved by Cramer's Rule to get $c_jx_j={w_j \over w_0}~(j=1,\ldots,m)$, where $w_0 \neq 0$ is the determinant of the matrix with entries ${D_ix_j \over x_j} ~(i=0,\ldots,m-1;~j=1,\ldots,m)$. Noting that this determinant is multilinear in the columns, we find that $h(w_0) \leq h(x_1;m-1)+\cdots+h(x_m;m-1)$. The same bound holds for the $h(w_j)~(j=1,\ldots,m)$. We conclude that $h(c_1x_1,\ldots,c_mx_m)=h({w_i \over w_0},\ldots,{w_m \over w_0})$ is at most
$$h(w_0)+h(w_1)+\cdots+h(w_m) ~\leq~ (m+1)\left(h(x_1;m-1)+\cdots+h(x_m;m-1)\right)$$
as required.
\bigskip
We deduce an explicit version of Lemma 4 of [Mass] (p.197), also without $G$.
\bigskip
\noindent
{\bf Lemma 3.2.} {\it For $m \geq 2$ suppose $x_0,x_1,\ldots x_m$ are in $K^*$ and linearly dependent over $C$ but $x_1,\ldots x_m$ are linearly independent over $C$. Then there is a relation
$$c_1x_1+ \cdots+c_mx_m~=~x_0 \eqno(3.2)$$
with $c_1,\ldots,c_m$ in $C$ and 
$$h\left({c_1x_1 \over x_0},\ldots,{c_mx_m \over x_0}\right) ~\leq~  (m+1)\left(h\left({x_1\over x_0};m-1\right)+\cdots+h\left({x_m\over x_0};m-1\right)\right).$$}
\bigskip
\noindent
{\it Proof.} There is certainly a relation (3.2) with $c_1,\ldots,c_m$ in $C$, and we apply Lemma 3.1 to the quotients ${x_1 \over x_0},\ldots,{x_m \over x_0}$. As $x_1,\ldots x_m$ are linearly independent over $C$, the conclusion (b) cannot hold. Now conclusion (a) is just what we need, and this completes the proof.
\bigskip
In section 5 we shall prove versions of Lemmas 3.1 and 3.2 that are uniform for $x_0,x_1,\ldots,x_m$ in  a finitely generated group $G$ as in [Mass]. By way of preparation, the next result illustrates the logarithmic nature of the quantities $h(~;i)$.
\bigskip
\noindent
{\bf Lemma 3.3.} {\it For any $x \neq 0,y \neq 0$ in $K$ and any integers $i \geq 0,e\geq 0$ we have $h(xy;i) \leq h(x;i)+h(y;i)$ and $h(x^e;i) \leq ih(x;i)$.}
\bigskip
\noindent
{\it Proof.} Let $D$ be in ${\cal D}(i)$. By distributing operators over the factors of $xy$ as in Leibniz, we see that ${D(xy) \over xy}$ is a sum with generalized binomial coefficients of products ${E(x) \over x}{F(y) \over y}$ with operators $E,F$ also in ${\cal D}(i)$. Taking $D=D_1,\ldots,D_l$ as in the definition of $h(xy;i)$, we deduce the first inequality of the present lemma by standard height calculations.
\medskip
When $e$ is a positive integer, a similar argument shows that ${D(x^e) \over x^e}$ is a sum with generalized binomial coefficients of products ${E_1(x) \over x}\cdots {E_e(x) \over x}$ with operators $E_1,\ldots,E_e$ also in ${\cal D}(i)$. Here $E_1\cdots E_e=D$, so that there are at most $i$ terms not equal to 1 in this product. Thus ${D(x^e) \over x^e}$ is a polynomial of total degree at most $i$ in the ${E(x) \over x}$ for $E$ in ${\cal D}(i)$. The second inequality now follows in a similar way, at least for $e \geq 1$. The result is trivial for $e=0$. 
\bigskip
\noindent
{\bf Lemma 3.4.} {\it For any $x \neq 0$ in $K$ and any integer $i \geq 0$ we have $h(x;i)\leq 4idh(x).$}
\bigskip
\noindent
{\it Proof.} This is trivial for $i=0$, so we assume from now on $i \geq 1$. We have an equation $A(x)=0$ as in Lemma 2.1, of degree $e \leq d$. Denote by $A'(t)$ the derivative with respect to $t$. Pick any $D$ in ${\cal D}(i)$. We claim that $B_i=(A'(x))^{2i-1}Dx$ is a polynomial in $x$ and various derivatives $D_0a$ of various coefficients $a$ of $A$, with coefficients in $k$ and of degree at most $(2i-1)(e-1)+1$ in $x$ and of total degree at most $2i-1$ in the $D_0a$. We prove this by induction on $i$.
\medskip
When $i=1$ we have for example $D={\partial \over \partial t_1}=\partial$ (say), and applying this to $A(x)=0$ yields $B_1=-\sum_{j=0}^e(\partial a_{e-j})x^j$ for which the claim is clear.
\medskip
Assuming $Dx={B_i \over (A'(x))^{2i-1}}$ with $B_i$ as above, we do the induction step by applying one more operator, again say ${\partial \over \partial t_1}=\partial$. We get
$$(A'(x))^{2i}\partial Dx~=~A'(x)\partial B_i - (2i-1)B_i\partial(A'(x)).$$
Here $\partial B_i$ involves $x$ to degree at most $(2i-1)(e-1)+1$ and also $x$ to degree at most $(2i-1)(e-1)$ multiplied by $\partial x={B_1 \over A'(x)}$, together with $D_0 a$ to total degree at most $2i-1$. Similarly $\partial(A'(x))$ involves $x$ to degree at most $e-1$ and also $x$ to degree at most $e-2$ (if $e \neq 1$) multiplied by $\partial x={B_1 \over A'(x)}$, together with $D_0 a$ to total degree at most $1$. Multiplying by $A'(x)$ we get $(A'(x))^{2i+1}\partial Dx$ involving $x$ to degree at most
$$e-1+\max\{(2i-1)(e-1)+1+(e-1),~(2i-1)(e-1)+e\}~=~(2(i+1)-1)(e-1)+1,$$
and the degree in $D_0a$ is at most $(2i-1)+1+1=2(i+1)-1$. This proves the claim in general.
\medskip
There follows at once the estimate
$$\log |B_i|_w ~\leq~ ((2i-1)(e-1)+1)\log\max \{1,|x|_w\}$$
for any $w$ not dividing $\infty$; and otherwise we get an extra term $(2i-1)\max\{\deg a_0,\ldots,\deg a_e\}$. The same estimates also hold for $\log |C|_w$ where $C=x(A'(x))^{2i-1}$. 
\medskip
Now write $B_{ij}$ for the $B_i$ corresponding to the operators $D_j~(j=1,\ldots,l)$ of ${\cal D}(i)$, so that ${D_jx \over x}={B_{ij} \over C}$. Then
$$h\left({D_1x \over x},\ldots,{D_lx \over x}\right)~=~\sum_w d_w \max\{\log |B_{i1}|_w,\ldots,\log |B_{il}|_w,\log |C|_w\}$$
which is at most
$$((2i-1)(e-1)+1)h(x)+(2i-1)d\max\{\deg a_0,\ldots,\deg a_e\}.$$
Finally by Lemma 2.1 this is at most
$$((2i-1)(e-1)+1)h(x)+(2i-1)eh(x) ~\leq~ 4ieh(x) ~\leq~ 4idh(x)$$
as required. This completes the proof of the present lemma.
\medskip
In view of our consistent use of the relative height (as opposed to the absolute height), the factor $d$ here looks like a normalization error. However it cannot be avoided, as the example $x=({t+1 \over t})^{1/d} ~(t=t_1)$ in $K=k(t)(x)=k(x)$ shows. One finds that the rational function ${1 \over x}{\partial^i x \over \partial t^i}$ has denominator $(t(t+1))^i$. So its height is at least $2id=2idh(x)$, which shows also that our dependence on $i$ is not too bad. Perhaps even the factor 4 essentially cannot be avoided.
\bigskip
\bigskip
\noindent
{\bf 4. Regulators.} Let $K$ be finitely generated and transcendental over $k$ as in the preceding section, and let $\cal B$ be a transcendence basis. Let $G$ be a subgroup of $K^*$ finitely generated modulo $k^*$; that is, $G/(G \cap k^*)$ is finitely generated. We show here how to define a regulator $R(G)=R_{\cal B}(G)$. 
\medskip
For all $w$ except finitely many we have $|g|_w = 1$ for every $g$ in $G$. Pick a set of $N \geq 1$ valuations containing these exceptions. We order the set to produce a map ${\cal L}$ from $G$ into ${\bf R}^N$ whose typical coordinate is $d_w \log |g|_w$. In fact by (2.2) ${\cal L}(G)$ lies in ${\bf Z}^N$ and is therefore discrete. Thus it is a (full) lattice in the real subspace it generates, whose dimension is the rank $r$ of $G/(G \cap k^*)$. If $r \geq 1$ we define the regulator just as the determinant 
$$R(G)=R_{\cal B}(G)=\det {\cal L}(G) \geq 1;$$
clearly independent of the set above or its ordering, and if $r=0$ we define $R(G)=1$. This does not quite coincide with the standard definition for the unit group in algebraic number theory, because the latter is obtained by a projection to one dimension lower. But they are equal up to a constant factor.
\medskip
The following example will be quoted later. With $K={\bf F}_p(t)$ (and the obvious $\cal B$) and $G_l$ generated by $t^l$ and $1-t$ we have $N=3$ corresponding to valuations at $t=0,1,\infty$; and so vectors $(l,0,l)$ and $(0,1,1)$ giving $R_{\cal B}(G_l)=l\sqrt{3}$.
\bigskip
\noindent
{\bf Lemma 4.1.} {\it Let $G,G'$ in $K^*$ be finitely generated modulo $k^*$ with $G$ of finite index in $G'$. Then
$$R(G)~=~{[G':G] \over [G' \cap k^*:G \cap k^*]}R(G')~=~[G'/(G'\cap k^*):G/(G \cap k^*)]R(G'),$$
where we identify $G/(G \cap k^*)$ as a subgroup of $G'/(G' \cap k^*)$.}
\bigskip
\noindent
{\it Proof}. The quotients $G/(G\cap k^*),G'/(G' \cap k^*)$ are torsion-free, both with the same rank, say $r$. If $r=0$ the lemma is trivial. Otherwise using elementary divisors we can find generators $\gamma_1,\ldots,\gamma_r$ of $G'/(G'\cap k^*)$ and positive integers $d_1,\ldots,d_r$ such that $\gamma_1^{d_1},\ldots,\gamma_r^{d_r}$ generate $G/(G \cap k^*)$. Then the relationship between ${\cal L}(G')$ and ${\cal L}(G)$ is clear, and the lemma follows.
\bigskip
\bigskip
\noindent
{\bf Lemma 4.2.} {\it Let $G$ in $K^*$ be finitely generated modulo $k^*$, let $x$ be in $K^*$, and let $G'$ be the group generated by $x$ and the elements of $G$. Then $R(G') \leq 2h(x)R(G)$.}
\bigskip
\noindent
{\it Proof}. It is geometrically clear that if $\Lambda$ is any lattice in euclidean space, then $\det(\Lambda+{\bf Z}{\bf v}) \leq \det(\Lambda)|{\bf v}|$ for the length, at least if $\bf v$ is not in the space spanned by $\Lambda$. But this continues to hold for all $\bf v$ provided only $|{\bf v}| \geq 1$ and $\Lambda+{\bf Z}{\bf v}$ remains discrete. In particular it holds for $\Lambda={\cal L}(G)$ and ${\bf v}={\cal L}(x)$. We conclude $R(G') \leq |{\cal L}(x)|R(G)$. Finally we have by definition and the product formula
$$h(x)~=~\sum_w\max\{0,m_w\}~=~{1 \over 2}\sum_w|m_w| \eqno(4.1)$$
for $m_w=d_w\log|x|_w$. And
$$|{\cal L}(x)|^2~=~\sum_wm_w^2 ~\leq~ (\sum_w|m_w|)^2~=~4(h(x))^2.$$
The lemma follows.
\bigskip
We can recover a basis from the regulator as follows. 
\bigskip
\noindent
{\bf Lemma 4.3.} {\it Let $G$ be a subgroup of $K^*$ finitely generated modulo $k^*$ with $G/(G \cap k^*)$ of rank $r \geq 1$. Then there are $g_1,\ldots,g_r$ in $G$ generating $G/(G \cap k^*)$, with
$$h(g_1)\cdots h(g_r) ~\leq~ {1 \over r}\delta(r)R(G)^2$$
for $\delta(r)=r^{3r}$.}
\bigskip
\noindent
{\it Proof}. By Minkowski's Second Theorem (see for example [Ca] Theorem V p.218) there are $\tilde g_1,\ldots,\tilde g_r$ in $G$ multiplicatively independent modulo $k^*$, with
$$|{\cal L}(\tilde g_1)|\cdots|{\cal L}(\tilde g_r)| ~\leq~{2^r \over V(r)}\det {\cal L}(G)~=~{2^r \over V(r)}R(G) \eqno(4.2)$$
for the Euclidean norms and the volume $V(r)$ of the unit ball in ${\bf R}^r$. By geometry $V(r) \geq ({2 \over \sqrt{r}})^r$. We get a basis in the standard way using the argument of Mahler-Weyl (see for example [Ca] Lemma 8 p.135); there results
$$|{\cal L}(g_i)| ~\leq~ \max\{1,{i \over 2}\}|{\cal L}(\tilde g_i)| ~~(i=1,\ldots,r),$$
and so ${2^r \over V(r)}$ in (4.2) gets replaced by ${r! \over 2^{r-1}}r^{r/2} \leq {r^{3r/2} \over 2^{r-1}}$. 
Now (4.1) gives
$$h(g)~=~\sum_w\max\{0,m_w\}~=~{1 \over 2}\sum_w|m_w|$$
for $m_w=d_w\log|g|_w$ in $\bf Z$. And $|m| \leq m^2$ for any $m$ in $\bf Z$, so we get
$$h(g) ~\leq~ {1 \over 2}\sum_wm_w^2~=~{1 \over 2}|{\cal L}(g)|^2.$$
Therefore 
$$h(g_1)\cdots h(g_r) ~\leq~ {4r^{3r} \over 2^{3r}}R(G)^2 ~<~ {1 \over r}\delta(r)R(G)^2$$ 
as desired.
\medskip
In view of (4.2) it seems a pity that the square of the regulator appears in Lemma 4.3. But it cannot be avoided. For example let $\alpha_1,\ldots,\alpha_l,\beta_1,\ldots,\beta_l$ be different constants in $k$, and consider $G$ generated by $g={(t-\alpha_1)\cdots(t-\alpha_l) \over (t-\beta_1)\cdots(t-\beta_l)}$ in $K=k(t)$. Then $R(G)=\sqrt{2l}$. The only possibilities for $g_1$ are $\gamma g^{\pm 1}$ with $\gamma$ constant. But then $h(g_1)=l$, so any bound $h(g_1) \leq \delta (1)R(G)$ is impossible. 
\bigskip
This leads to the following uniform version of Lemma 3.4 when $x$ lies in $G$. Write $G_k$ for the group generated by the elements of $G$ and $k^*$. 
\bigskip
\noindent
{\bf Lemma 4.4.} {\it Let $G$ be a subgroup of $K^*$ finitely generated modulo $k^*$ with $G/(G \cap k^*)$ of rank $r \geq 1$. Then for any $g$ in $G$ we have $h(g;i) \leq 4i^2d\delta(r)R(G)^2.$ Further for any positive integer $l$ there is $g_0$ in $G_k$ and $g'$ in $G$ with $g=g_0g'^l$ and $h(g_0) \leq l\delta(r)R(G)^2.$}
\bigskip
\noindent
{\it Proof.} Choose basis elements $g_1,\ldots,g_r$ according to Lemma 4.3, and write $g=cg_1^{e_1}\cdots g_r^{e_r}$ for rational integers $e_1,\ldots,e_r$ and $c$ in $k^*$. Replacing some of the $g_j$ by their inverses, we can assume that all $e_j \geq 0$; this does not affect the estimate in Lemma 4.3. Then by Lemma 3.3
$$h(g;i)~=~h(g_1^{e_1}\cdots g_r^{e_r};i) ~\leq~ h(g_1^{e_1};i)+\cdots+h(g_r^{e_r};i) ~\leq~ i(h(g_1;i)+\cdots+h(g_r;i)).$$
This in turn by Lemma 3.4 is at most 
$$4i^2d(h(g_1)+\cdots +h(g_r)) ~\leq~ 4i^2drh(g_1)\cdots h(g_r) ~\leq~4i^2d\delta(r)R(G)^2 \eqno(4.3)$$
as required in the first part of the present lemma. And the second part follows by writing $e_j=f_j+le_j'$ with $0 \leq f_j < l~(j=1,\ldots,r)$ (compare also [D] p.197), taking $g_0=cg_1^{f_1}\cdots g_r^{f_r}, g'=g_1^{e_1'}\cdots g_r^{e_r'}$ and using the inequality in (4.3).
\bigskip
The final result of this section will lead easily to a quantitative version of Lemma 2 of [Mass] (p.193), such as those mentioned in [Mass] (pp 194,195). However it involves better constants and is no longer restricted to positive characteristic. It is here, by the way, that the radical $\sqrt{G}$ makes its essential appearance in the whole story. 
\bigskip
\noindent
{\bf Lemma 4.5.} {\it Suppose that $x,y$ are in $K^*$ with $x$ not in $\sqrt{G_k}$ and ${y^q \over x}$ in $G$ for some positive integer $q$. Then $q \leq 2h(x)R(G)$.}
\bigskip
\noindent
{\it Proof.} Let $G'$ be the group generated by $x$ and the elements of $G$, and let $G''$ be the group generated by $y$ and the elements of $G$, so that $G'$ lies in $G''$. Since $x$ is not in $\sqrt{G}$, it is easy to see that the index $[G'':G']=q$. Since $x$ is not even in $\sqrt{G_k}$, it is even easier to see that 
$G \cap k^*=G' \cap k^*=G'' \cap k^*$. Thus by Lemma 4.1 we have $R(G')=qR(G'') \geq q$. On the other hand $R(G') \leq 2h(x)R(G)$ by Lemma 4.2, and the result follows.
\bigskip 
\bigskip
\noindent
{\bf 5. Dependence with regulators.} Let $K$ be finitely generated and transcendental over $k$ as in the preceding sections, and let $\cal B$ be a transcendence basis, now assumed separable, with elements $t_1,\ldots,t_b$. We continue to abbreviate the height $h_{\cal B}$ as $h$, and again we write $C$ for the field of differential constants of $K$.  
\medskip
The following result eliminates the height functions $h(x,m-1)$ from Lemma 3.1, thereby providing a more useful explicit version of Lemma 3 of [Mass].
\bigskip
\noindent
{\bf Lemma 5.1.} {\it Let $G$ in $K^*$ be finitely generated of rank $r \geq 1$ modulo $k^*$, and for $m \geq 2$ suppose $c_1,\ldots,c_m$ are in $C$ and $g_1,\ldots g_m$ are in $G$ with 
$$c_1g_1+ \cdots+c_mg_m~=~1.$$
Then either

(a) $h(c_1g_1,\ldots,c_mg_m) ~\leq~ 4m^4d\delta(r)R(G)^2$

\noindent
or 

(b) $g_1,\ldots,g_m$ are linearly dependent over $C$.}
\bigskip
\noindent
{\it Proof.} Just use Lemma 3.1 together with the inequalities 
$$h(g;m-1) ~\leq~ 4(m-1)^2d\delta(r)R(G)^2\eqno(5.1)$$ 
from Lemma 4.4, with $g=g_1,\ldots,g_m$. 
\bigskip
Similarly we deduce a more useful explicit version of Lemma 4 of [Mass].
\bigskip
\noindent
{\bf Lemma 5.2.} {\it Let $G$ in $K^*$ be finitely generated of rank $r \geq 1$ modulo $k^*$, and for $m \geq 2$ suppose $g_0,g_1,\ldots g_m$ are in $G$ and linearly dependent over $C$ but $g_1,\ldots g_m$ are linearly independent over $C$. Then there is a relation
$$c_1g_1+ \cdots+c_mg_m=g_0$$
with $c_1,\ldots,c_m$ in $C$ and 
$$h\left({c_1g_1 \over g_0},\ldots,{c_mg_m \over g_0}\right) ~\leq~ 4m^4d\delta(r)R(G)^2.$$}
\bigskip
\noindent
{\it Proof.} Just use Lemma 3.2 and (5.1), this time with $g={g_1 \over g_0},\ldots,{g_m \over g_0}$. 
\bigskip
We have followed the proof in [Mass] quite closely. It would have been nice to see the well-known number ${m(m-1) \over 2}$ in place of $4m^4$, and also some notion of genus and $S$-units as in various formulations of $abc$ matters over function fields. But despite the considerations of Chapter 14 of [BG] in zero characteristic and those of Hsia and Wang [HW] for arbitrary characteristic we have been unable to supply a satisfactory version. The results of [HW] are especially interesting in their use of divided derivatives or hyperderivations, which for example in characteristic $p$ leads to linear dependence over the field of $p^e$th powers, not just over $C$ with $e=1$. If this could be done in our situation, then it would probably lead to simplifications in the rest of our proof, and possibly to the elimination of the Proposition in section 8. But it seems that the results of [HW] cannot be directly applied to our Lemma 5.1, due to the presence of $c_1,\ldots,c_m$ whose heights cannot be controlled.
\bigskip 
\bigskip
\noindent
{\bf 6. Isotriviality.} We take a well-earned break from estimating. From now on $K$ will have positive characteristic $p$ (actually this assumption is not really needed until section 8), and, as in section 1, we write ${\bf F}_K$ for $\overline {{\bf F}_p} \cap K$. This field plays the role of $k$ in sections 2,3,4,5. 
\medskip
Suppose $n \geq m \geq 1$. For $a(i,j)$ in $K$ the normalized equations 
$$X_i~=~a(i,0)X_0+\cdots+a(i,m-1)X_{m-1}~=~\sum_{j=0}^{m-1}a(i,j)X_j ~~(i=m,m+1,\ldots,n) \eqno(6.1)$$
define in ${\bf P}_n$ a linear variety $V$ of dimension $m-1$. When $G$ is a subgroup of $K^*$, we need some conditions which ensure that $V$ is $G$-isotrivial. 
\medskip
Now any $G$-automorphism taking $(X_0, \ldots ,X_n)$ to $(g_0X_0, \ldots ,g_nX_n)$ leads after renormalization to new coefficients ${g_i \over g_j}a(i,j)$. If the new forms are defined over ${\bf F}_K$, then every non-zero $a(i,j)$ has the shape ${g_j \over g_i}\alpha(i,j)$ for non-zero $\alpha(i,j)$ in ${\bf F}_K$. In particular each equation in (6.1) defines a $G$-isotrivial variety. But also each quotient 
$${a(i_1,j_1)a(i_2,j_2)a(i_3,j_3) \cdots a(i_{k-1},j_{k-1})a(i_k,j_k) \over a(i_1,j_2)a(i_2,j_3)a(i_3,j_4) ~~\cdots~~ a(i_{k-1},j_k)a(i_k,j_1)} ~~~(k=2,\ldots,n+1),\eqno(6.2)$$
with everything in the numerator and denominator non-zero, lies in ${\bf F}_K$. The following result gives a converse statement which guarantees that the equations (6.1) become defined over ${\bf F}_K$ after applying a suitable $G$-automorphism and renormalizing. In particular it guarantees that $V$ is $G$-isotrivial; but without the need to recombine the equations.
\bigskip
\noindent
{\bf Lemma 6.1.} {\it Suppose that each equation in (6.1) defines a $G$-isotrivial variety, and that each quotient (6.2) lies in ${\bf F}_K$ provided everything in the numerator and denominator is non-zero. Then $V$ is $G$-isotrivial.}
\bigskip
\noindent
{\it Proof.} We argue by induction on the number $n-m+1 \geq 1$ of equations. If $n-m+1=1$ then the result is obvious without using (6.2). Suppose we have done it for $n-m \geq 1$ equations,  namely the first $n-m$ in (6.1), and let us add another equation, namely the last one in (6.1). 
\medskip
Restricting to $i < n$ and the appropriate $j$ in (6.2), we see from the induction hypothesis that a suitable $G$-automorphism trivializes the first $n-m$ equations, without bothering about $X_n$. This means that we can assume that all $a(i,j) \neq 0 ~(i<n)$ are in ${\bf F}_K$; while the isotriviality of the last equation means that all $a(n,j) \neq 0$ are in $G$. We now want to trivialize the last equation.
\medskip
How can we trivialize a given coefficient $a(n,j) \neq 0$ in the last equation? If all $a(i,j)=0 ~(i<n)$, so that the first $n-m$ equations did not involve $X_j$, then we can simply replace $X_j$ by $a(n,j)X_j$  and this will not change the first $n-m$ equations. We do this for all such $j$. 
\medskip
If there is only a single $j$ with some $a(i,j) \neq 0 ~(i<n)$, then we can still replace $X_j$ by $a(n,j)X_j$; but we then have to correct the new coefficients ${a(i,j) \over a(n,j)} \neq 0$ of $X_j$ in the $i$th equation by replacing $X_i$ by $a(n,j)X_i~(i=m,\ldots,n-1)$. Things are less easy when there is more than one such $j$. Call these ``bad".
\medskip
Now we say for different $j,j'$ in the set $\{0, \ldots ,m-1\}$ that $j \sim j'$ if there is $i<n$ with 
$$a(i,j)a(i,j') \neq 0 \eqno(6.3)$$ 
(in particular then $j,j'$ are both bad). This relation is symmetric but probably not transitive. We can extend it via reflexivity and transitivity to a genuine equivalence relation on the bad elements of $\{0, \ldots ,m-1\}$, which we then denote 
by $\simeq$.
\medskip
We assume for the moment that there is a single equivalence class: any two $j,j'$ are related.
\medskip
Let $j,j'$ be different bad elements, so that $a(i,j) \neq 0,~ a(i',j') \neq 0$ for some $i,i' <n$. From our equivalence class assumption $j \simeq j'$. Suppose that 
$$j=j_1 \sim j_2  \sim \cdots  \sim j_{k-1} \sim j_k = j',$$
where of course we can take $2 \leq k \leq n+1$. Then we get from (6.3)
$$a(i_1,j_1)a(i_1,j_2) \neq 0, ~a(i_2,j_2)a(i_2,j_3) \neq 0, \ldots, ~a(i_{k-1},j_{k-1})a(i_{k-1},j_k) \neq 0 $$
for some $i_1,i_2,\ldots,i_{k-1} <n$. We use (6.2) with $i_k=n$ to see that
$${a(i_1,j_1)a(i_2,j_2)a(i_3,j_3) \cdots a(i_{k-1},j_{k-1})a(n,j') \over a(i_1,j_2)a(i_2,j_3)a(i_3,j_4) ~~\cdots~~ a(i_{k-1},j_k)a(n,j)}$$
lies in ${\bf F}_K$. However the first $k-1$ terms in both numerator and denominator already lie in ${\bf F}_K$, because we trivialized the first $n-m$ equations. Consequently ${a(n,j') \over a(n,j)}$ lies in ${\bf F}_K$.
\medskip
Thus we have shown that all $a(n,j)$ for bad $j$ are multiples of a single one, call it $g$, by elements of ${\bf F}_K$. Now they can be simultaneously trivialized on replacing $X_j$ by $gX_j$. Again we must correct the new coefficients ${a(i,j) \over g} \neq 0$ of $X_j$ in the $i$th equation by replacing $X_i$ by $gX_i~(i=m,\ldots,n-1)$.
\medskip
What happens if there is more than a single equivalence class on the bad elements of $\{0, \ldots ,m-1\}$? Say there are $h \geq 2$ classes $J_1 \ldots,J_h$. Let $I_1$ be the set of $i$ in $\{m,\ldots,n-1\}$ for which there is $j$ in $J_1$ with $a(i,j) \neq 0$; and similarly for $I_2,\ldots,I_h$. Then $I_1,I_2,\ldots,I_h$ are disjoint, because for example with any $j_1$ in $J_1$ and any $j_2$ in $J_2$ there can be no $i$ with $a(i,j_1)a(i,j_2) \neq 0$, else by (6.3) we would have $j_1 \sim j_2$. (If one wishes, one can convert the matrix of the first $n-m$ equations into a block matrix using row and column permutations.) The argument above, using $i_1,\ldots,i_{k-1}$ in $I_1$, shows that all non-zero $a(n,j)~ (j \in J_1)$ are multiples of a single one, call it $g_1$, by elements of ${\bf F}_K$. Similarly we get $g_2,\ldots,g_h$. Now we can trivialize the last row as follows. We replace the $X_j ~(j \in J_1)$ by $g_1X_j $ and we correct the effect by replacing $X_i$ by $g_1X_i~(i \in I_1)$. Similarly using $g_2,\ldots,g_h$ we trivialize the remaining coefficients. This completes the proof.
\bigskip 
\bigskip
\noindent
{\bf 7. Automorphisms.} As above let $K$ be a field, finitely generated and transcendental over ${\bf F}_p$, with $G$ a subgroup of $K^*$. Suppose a linear variety in ${\bf P}_n$ is defined over $K$ and $G$-isotrivial. Then by definition there is a $G$-automorphism $\psi$ taking it to something defined over ${\bf F}_K=\overline{{\bf F}_p} \cap K$. To make our Theorems 1,2 and 3 fully effective we have to estimate this $\psi$; indeed after doing the whole descent to single points using Theorem 1, for example, it is mainly $G$-automorphisms that are left.
\medskip
Now it is convenient to use the projective height $h^{\bf P}=h_{\cal B}^{\bf P}$ defined on ${\bf P}_{l-1}(K)$ by
$$h^{\bf P}(x_1,\ldots,x_l)~=~\log \prod_w \max\{|x_1|_w^{d_w},\ldots,|x_l|_w^{d_w}\}.$$
This yields at once a height $h(\psi)$ of a $G$-automorphism $\psi$, defined by (1.7), as
$$h(\psi)=h^{\bf P}(g_0,\ldots,g_n).$$
Also if $V$ is linear in ${\bf P}_n$ defined over $K$, it yields a height $h(V)$ in the standard way via the Grassmannian coordinates of $V$; see for example [S] (p.28), which however is in the context of number fields with euclidean norms at the archimedean valuations. Here we have no archimedean valuations, so the norm problem is irrelevant. If $m-1 \geq 0$ is the dimension of $V$, then its Grassmannians $A(I)$ correspond to subsets $I$ of $\{0,\ldots,n\}$ with cardinality $n-m+1 \leq n$. The Northcott Property extends at once to this height. Also for $\psi$ in (1.7) the Grassmannians of $\psi(V)$ are the ${A(I) \over g(I)}$, where $g(I)=\prod_{i \in I}g_i$. It follows easily that
$$h(\psi(V))~\leq~ h(V)+nh(\psi),~~h(\psi^{-1})~\leq~ nh(\psi). \eqno(7.1)$$
Less obvious is the following, which involves a second linear variety $W$ also over $K$.
\bigskip
\noindent
{\bf Lemma 7.1.} {\it If $V \cap W$ is non-empty then we have $h(V \cap W) \leq h(V)+h(W)$. If further $X_{n-1}\neq 0$ on $V$ and the equations of $V$ do not involve $X_n$, and $W$ is defined by $X_n=aX_{n-1}$ then $h(V \cap W) \geq \max\{h(V),h(W)\}$.}
\bigskip
\noindent
{\it Proof.} The upper bound may be compared with the inequality $h(V \cap W)+h(V+W) \leq h(V)+h(W)$ due independently to Struppeck-Vaaler [SV] (Theorem 1 p.493) and Schmidt [S] (Lemma 8A p.28). These are proved over number fields; however it is easily checked that the proof in [S] remains valid with trivial modifications. Already a special case was noted by Thunder [T] whose Lemma 5 (p.157) implies $h(V+W) \leq h(V)+h(W)$ over function fields of a single variable provided $V \cap W$ is empty.
\medskip
Regarding the lower bound, let $A(I)$ be the Grassmannians of $V$. Then it is easy to verify that the Grassmannians of $V \cap W$ consist of the $A(I)$ together with the $aA(J)$ for $J$ not containing $n-1$. There follows $h(V \cap W) \geq h(V)$ at once. Also $X_{n-1}\neq 0$ on $V$ means that at least one $A=A(J)$ is non-zero (see for example Theorem 1 of [HP] p.298), so we get also $h(V \cap W) \geq h^{\bf P}(A,aA)=h(a)=h(W)$. This completes the proof.
\bigskip
It is the following result which enables $\psi$ to be estimated in the Descent Steps.
\bigskip
\noindent
{\bf Lemma 7.2.} {\it Suppose that $V$ is defined over $K$ and is $G$-isotrivial. Then there is a $G$-automorphism $\psi$ with $\psi(V)$ defined over ${\bf F}_K$ and $h(\psi) \leq n!h(V)$.}
\bigskip
\noindent
{\it Proof.} Suppose $\dim V=m-1$ with Grassmannians $A(I)$; then as noted above the Grassmannians of $\psi(V)$ are the ${A(I) \over g(I)}$, where $g(I)=\prod_{i \in I}g_i$. If $\psi(V)$ is defined over ${\bf F}_K$ then these have the shape $\lambda \alpha(I)$ for $\lambda$ in $K^*$ and $\alpha(I)$ in ${\bf F}_K$. Thus we have $A(I)=\lambda \alpha(I)g(I)$ for all $I$; but we can restrict to the set $\cal I$ of all $I$ with $A(I) \neq 0$ (and so $\alpha(I) \neq 0$). We can eliminate the $\lambda$ by fixing $I_0$ in $\cal I$; this gives
$${g(I) \over g(I_0)}~=~{A(I) \over A(I_0)}{\alpha(I_0) \over \alpha(I)}~~~(I \in {\cal I}). \eqno(7.2)$$
Conversely (7.2) implies that $\psi(V)$ is defined over ${\bf F}_K$.
\medskip
To solve (7.2) for $g_0,\ldots,g_n$ we divide the numerator and denominator of the left-hand side by $g_0^{n-m+1}$ and write it as $({g_1 \over g_0})^{a(I,1)}\cdots ({g_n \over g_0})^{a(I,n)}$ for integers $a(I,i)$ which are $0,1,-1$. If the vectors ${\bf a}(I)~(I \in {\cal I})$ with coordinates $a(I,i)~(i=1\ldots,n)$ have full rank $n$ then we can solve (7.2) by choosing ${\bf a}(I_1),\ldots,{\bf a}(I_n)$ linearly independent and then solving the subset of (7.2) with $I=I_1,\ldots,I_n$. A multiplicative form of Cramer's Rule gives
$$\left({g_i \over g_0}\right)^b=Q_1^{b_{i1}}\cdots Q_n^{b_{in}},~~~Q_j={A(I_j) \over A(I_0)}{\alpha(I_0) \over \alpha(I_j)}~~~(j=1,\ldots,n)$$
with integers $b \neq 0$ and $b_{ij}$. These $b_{ij}$ are minors of a matrix with entries $0,1,-1$ and so $|b_{ij}| \leq (n-1)!$.
\medskip
Now taking heights leads to
$$|b|h\left({g_1 \over g_0},\ldots,{g_n \over g_0}\right) ~\leq~ \max_{i=1,\ldots, n}\{|b_{i1}|+\cdots+|b_{in}|\}h(Q_1,\ldots,Q_n).$$
The height on the left is $h(\psi)$ and that on the right at most $h(V)$. The result follows at once, at least under our assumption that the ${\bf a}(I)~(I \in {\cal I})$ have full rank $n$.
\medskip
If this assumption does not hold, then we simply increase the rank by successively adjoining unit vectors ${\bf e}_k$ until the rank becomes $n$; this amounts to the addition of equations ${g_k \over g_0}=1$. Now we take a subset of $n$ independent equations and solve again with Cramer. The resulting estimates are certainly no larger than before, and this completes the proof.
\bigskip 
\bigskip
\noindent
{\bf 8. A proposition.} This, the main result of the present section, is a first step in the proof of the Descent Step over $\sqrt{G}$, with $V$ in ${\bf P}_n~(n \geq 2)$ either a hyperplane or defined over a finite field. We continue with our assumption that $K$ is finitely generated over ${\bf F}_p$; thus ${\bf F}_K=\overline{{\bf F}_p} \cap K$ is a finite field. Let $G$ in $K^*$ be finitely generated of rank $r \geq 1$ modulo ${\bf F}_K^*$; now we may write without confusion simply that $G$ is finitely generated. It is known that the radical $\sqrt{G}$, which by definition lies still in $K$, is also finitely generated (see for example [Mass] p.195), also clearly of rank $r$ over ${\bf F}_K^*$. For the moment we work exclusively with this radical. We further assume that $K$ is transcendental over ${\bf F}_p$ and we choose any separable transcendence basis $\cal B$; then we are free to apply the results of sections 3,4 and 5 about heights $h=h_{\cal B}$ and regulators $R=R_{\cal B}$.
\medskip
We say that $V$ is transversal if every coordinate $X_i ~(i=0,\ldots,n)$ actually occurs in the defining equations. This property is independent of the choice of equations. Its purpose is to prevent ``free variables" as in (1.1) with $a_i \neq 0$.  
\medskip
Transversality is a harmless restriction because we could overcome it simply by working in lower dimensions. Clearly every linear subvariety of a transversal variety is also transversal. Also a transversal variety must be proper (i.e. not the full ${\bf P}_n$). 
\medskip
We recall the function $\delta$ from Lemma 4.3.
\bigskip
\noindent
{\bf Proposition}. {\it Let $V$ be a transversal linear subvariety of ${\bf P}_n$ defined over $K$, and suppose either that $V$ has dimension $n-1$ or that $V$ is defined over some ${\bf F}_q$. Suppose also that $V$ is not contained in any coset $T \neq {\bf P}_n$. Let $\pi$ be any point of $V(\sqrt{G})$. 

If $V$ has dimension $n-1$, then either

\noindent
(i) there is a proper linear subvariety $W$ of $V$, also defined over $K$, with 
$$h(W) ~\leq~ 8n^54^nd\delta(n+r)h(V)^{2n}R(\sqrt{G})^2,$$
such that $\pi$ lies in $W(\sqrt{G})$,

\noindent
or

\noindent
(ii) there is a $\sqrt{G}$-automorphism $\psi$ with 
$$h(\psi) ~\leq~ np\delta(n+r)R(\sqrt{G})^2,$$
a point $\pi'$ and a linear subvariety $V'$ of ${\bf P}_n$ such that $\pi=\psi(\pi'^p)$ and $V=\psi(V'^p)$.

If $V$ is defined over ${\bf F}_q$, then either 

\noindent
(i) there is a proper linear subvariety $W$ of $V$, also defined over $K$, with 
$$h(W) ~\leq~ 8n^54^nd\delta(n+r)R(\sqrt{G})^2,$$
such that $\pi$ lies in $W(\sqrt{G})$,

\noindent
or

\noindent
(iii) there is a point $\pi'$ in ${\bf P}_n(\sqrt{G})$ with $\pi=\pi'^p$.}
\bigskip
\noindent
{\it Proof.} Suppose first that $V$ has dimension $n-1$. Then we just have to follow the arguments of the proof of Lemma 5 of [Mass] (p.197). Because these arguments are expressed in terms of ``broad sets" and this notion is no longer appropriate, we write out all the details.
\medskip
Because $V$ is transversal, we may work affinely with a point $\pi=(x_1,\ldots,x_n)$ satisfying a single equation
$$a_1x_1+\cdots+a_nx_n=1 \eqno(8.1)$$
with non-zero coefficients. As in section 3 write $C$ for the field of $p$th powers in $K$, and consider
$$s=\dim_C(Ca_1x_1+\cdots+Ca_nx_n),$$
so that $1 \leq s \leq n$.
\medskip
First suppose that $s=n$. Then we apply Lemma 5.1 with $k={\bf F}_K$, $m=n$ and $c_1=\cdots=c_m=1$ and $g_1=a_1x_1,\ldots,g_m=a_mx_m$. So the group must be enlarged by adjoining $a_1,\ldots,a_n$ to $\sqrt{G}$, becoming of rank at most $n+r$. The enlarged regulator $R$ can be estimated by Lemma 4.2, and we find
$$R ~\leq~ 2^nh(a_1)\cdots h(a_n)R(\sqrt{G}) ~\leq~ 2^nh(V)^nR(\sqrt{G}). \eqno(8.2)$$
The conclusion $(b)$ of Lemma 5.1 is ruled out by $s=n$; and the conclusion $(a)$ shows that 
$$h(a_1x_1,\ldots,a_nx_n)~\leq~ 4n^4d\delta(n+r)R^2.$$ 
It follows that $h(\pi)=h(x_1,\ldots,x_n)$ is at most
$$4n^4d\delta(n+r)R^2+h(a_1^{-1},\ldots,a_n^{-1})~\leq~ 4n^4d\delta(n+r)R^2+nh(V)$$
and so from (8.2) we deduce
$$h(\pi)~\leq~4n^44^nd\delta(n+r)h(V)^{2n}R(\sqrt{G})^2+nh(V)~\leq~8n^44^nd\delta(n+r)h(V)^{2n}R(\sqrt{G})^2. \eqno(8.3)$$ 
So this gives $W=\{\pi\}$ for $(i)$ of the Proposition; and for these $h(W)=h(\pi)$ is bounded as in (8.3).
\medskip
Next suppose that $1 <s <n$. By means of a permutation we can assume that $g_1=a_1x_1,\ldots,g_s=a_sx_s$ are linearly independent over $C$. Take any $k$ with $s+1 \leq k \leq n$; then we can apply Lemma 5.2 with $m=s$ and $g_0=a_kx_k$, $\sqrt{G}$ being enlarged as above. We find relations
$$\sum_{j=1}^sc_{kj}a_jx_j~=~a_kx_k ~~~(k=s+1,\ldots,n) \eqno(8.4)$$
with $c_{kj}$ in $C$ and the quotients
$$f_{kj}~=~c_{kj}{a_jx_j \over a_kx_k} ~~~(j=1,\ldots,s;~k=s+1,\ldots,n) \eqno(8.5)$$
satisfying
$$h(f_{k1},\ldots,f_{ks}) ~\leq~  4s^4d\delta(n+r)R^2~~~(k=s+1,\ldots,n) \eqno(8.6)$$
We use (8.4) to eliminate the $a_kx_k ~(k=s+1,\ldots, n)$ in (8.1). We find
$$c_1a_1x_1+ \cdots+c_sa_sx_s~=~1 \eqno(8.7)$$
with 
$$c_j~=~1+\sum_{k=s+1}^nc_{kj} ~~~(j=1,\ldots,s) \eqno(8.8)$$
also in $C$.
\medskip
Next apply Lemma 5.1 with $m=s$ to (8.7) and $g_j=a_jx_j~(j=1,\ldots,s)$ also in the enlarged $\sqrt{G}$. Again conclusion $(b)$ is impossible. It follows that the
$$f_j~=~c_ja_jx_j~~~(j=1,\ldots,s) \eqno(8.9)$$
satisfy
$$h(f_1,\ldots,f_s) ~\leq~ 4s^4d\delta(n+r)R^2.\eqno(8.10)$$
\medskip
So in (8.5) certain quotients ${x_j \over x_k}$ are bounded modulo $C$ whereas in (8.9) certain $x_j$ themselves are bounded modulo $C$. We can eliminate $C$ by substituting (8.8) into (8.9) and using (8.5) to get
$$f_j~=~a_jx_j+\sum_{k=s+1}^nf_{kj}a_kx_k ~~(j=1,\ldots,s). \eqno(8.11)$$
Since $a_j \neq 0~(j=1,\ldots,s)$ these express the fact that $\pi=(x_1,\ldots,x_n)$ lies on a linear variety $V'$ of dimension $n-s$; and because $s \neq 1$ this dimension is strictly less than the dimension $n-1$ of $V$. So we can take $W$ as the intersection of $V'$ with $V$. This is in fact $V'$ because if we add up all the above equations $(8.11)$ and use (8.4),(8.5),(8.7),(8.9), then we end up with (8.1). 
\medskip
Now we have to estimate the height of (8.11). In the corresponding matrix, every column has by (8.6) and (8.10) height at most $4s^4d\delta(n+r)R^2+h(V)$, which as above in (8.3) we can estimate by $B=8n^44^nd\delta(n+r)h(V)^{2n}R(\sqrt{G})^2$. It follows that
$$h(W) ~\leq~ sB ~\leq~ 8n^54^nd\delta(n+r)h(V)^{2n}R(\sqrt{G})^2.$$
This too settles $(i)$ of the Proposition.
\medskip
Finally suppose $s=1$. This means that $a_1x_1,\ldots,a_nx_n$ are in $C$. By Lemma 4.4 with $l=p$ we can write $x_j=g_jx_j'^p$ with $g_j,x_j'$ in $\sqrt{G} ~(j=1,\ldots,n)$ and 
$$h(g_j) ~\leq~ p\delta(r)R(\sqrt{G})^2 ~\leq~ p\delta(n+r)R(\sqrt{G})^2~~(j=1,\ldots,n).$$  
Then $a_jg_j$ is in $C$ so has the form $a_j'^p ~(j=1,\ldots,n)$. Finally
$$1~=~a_1x_1+\cdots+a_nx_n~=~a_1'^px_1'^p+\cdots+a_n'^px_n'^p~=~(a_1'x_1'+\cdots+a_n'x_n')^p,$$ 
and this gives part $(ii)$ of the Proposition, with $\psi$ as in (1.7) above for $g_0=1$, $\pi'=(x_1',\ldots,x_n')$, and $V'$ defined by (8.1) above with the new coefficients $a_1',\ldots,a_n'$.
\medskip
This proves the Proposition when $V$ has dimension $n-1$. Incidentally when the coefficients in (8.1) are in some ${\bf F}_q$, then the argument for $s=1$ shows that $x_1,\ldots,x_n$ are in $C$. So they are $p$-th powers $x_1'^p,\ldots,x_n'^p$; and clearly $x_1',\ldots,x_n'$ are in $\sqrt{G}$. Thus we get the conclusion $(iii)$ of the Proposition when $V$ has dimension $n-1$. And the case $s \neq 1$ leads of course to $(i)$. So it remains only to treat $V$ of dimension $m-1 < n-1$ defined over some ${\bf F}_q$. 
\medskip
This we do by expressing the affine equations of $V$ in triangular form, which after a permutation we can suppose are
$$x_i~=~a_{i0}+a_{i1}x_1+\cdots+a_{i,m-1}x_{m-1}~~(i=m,m+1,\ldots,n) \eqno(8.12)$$
with the $a_{ij}$ in ${\bf F}_q$. This gives $V=V_m \cap \cdots \cap V_n$ for the varieties defined individually by each equation.
\medskip
Consider the first equation. There may be some zero coefficients $a_{mj}$, but not all are zero, because $V(\sqrt{G})$ is non-empty. In fact at least two are non-zero otherwise $V$would be contained in a coset $T \neq {\bf P}_n$ contrary to our assumption. We can thus regard $V_m$ as a transversal variety of codimension 1 in some projective space of dimension at least 2 and at most $m<n$. Applying the Proposition for the cases already proved, we get two possibilities $(i),(iii)$. If $(i)$ holds for $V_m$, then we get a proper subvariety $W_m$ of $V_m$ with 
$$h(W_m) ~\leq~ 8n^54^nd\delta(n+r)R(\sqrt{G})^2. \eqno(8.13)$$
But it is not difficult to see that each $W_m$ intersects the remaining intersection $U_m=\bigcap_{i \neq m}V_{i}$ in a proper subspace of $V=V_m \cap U_m$. For example the triangular nature of (8.12) makes it clear that $x_{m+1},\ldots,x_n$ are determined by $x_1,\ldots,x_{m-1}$ on $U_m$, and then that $x_m$ is determined by $x_1,\ldots,x_{m-1}$ on $W_m$ in $V_m$; but also some non-zero polynomial of degree at most 1 in $x_1,\ldots,x_{m-1}$ must vanish on $W_m$. So $W=W_m \cap U_m$ has dimension strictly less than $m-1$. By Lemma 7.1 we have $h(W) \leq h(W_m)$. So by (8.13) we get $(i)$ of the Proposition for the original $V$. But what happens if $(iii)$ holds for $V_m$?
\medskip
This means that all the $x_j$ actually occurring in the first equation of (8.12) are $p$-th powers, which certainly goes some way in the direction of $(iii)$ for $V$. But then we can try the second equation instead. Either we get a $W$ as above, or all the $x_j$ actually occurring in the second equation of (8.12) are $p$-th powers. And so on. In the end, we either get $W$ or that all the $x_j$ actually occurring in all the equations (8.12) are $p$-th powers. Because $V$ is transversal this does give the full $(iii)$ for $V$; and so completes the proof of the Proposition. 
\bigskip 
\bigskip
\noindent
{\bf 9. The main estimate.} This is a quantitative version of our Descent Step over $\sqrt{G}$ without the requirement that the subvarieties $W$ are isotrivial. This leads to a relatively small exponent attached to the height $h(V)$. As before $n \geq 2$, and we continue with our assumption that $K$ is finitely generated and transcendental over ${\bf F}_p$, with separable transcendence basis $\cal B$ and ${\bf F}_K=\overline{{\bf F}_p} \cap K$; further $G$ is finitely generated of rank $r \geq 1$ modulo ${\bf F}_K^*$.
\bigskip
\noindent
{\bf Main Estimate}. {\it Let $V$ be a positive-dimensional linear subvariety of ${\bf P}_n$ defined over $K$ but not a coset. 

\noindent
(a) If $V$ is not $\sqrt{G}$-isotrivial, then
$$V(\sqrt{G})~=\bigcup_{W \in {\cal W}}W(\sqrt{G})$$
for a finite set $\cal W$ of proper linear subvarieties $W$ of $V$, also defined over $K$ and with
$$h(W) ~\leq~  8n^2d(10n^3\delta(n+r))^{2n+1}h(V)^{2n}R(\sqrt{G})^{6n+2}.$$

\noindent
(b) If $V$ is $\sqrt{G}$-isotrivial and $\psi(V)$ is defined over ${\bf F}_q$, then
$$V(\sqrt{G})~=~\psi^{-1}\left(\bigcup_{W \in {\cal W}}~\bigcup_{e=0}^\infty(\psi(W)(\sqrt{G}))^{q^e}\right)$$
for a finite set $\cal W$ of proper linear subvarieties $W$ of $V$, also defined over $K$ and with
$$h(\psi(W)) ~\leq~  8n^54^n(q/p)d\delta(n+r)R(\sqrt{G})^2.$$}
\bigskip
\noindent
{\it Proof.} We prove this first when $V$ is transversal and not contained in any coset $T \neq {\bf P}_n$. 
\medskip
We start with $\sqrt{G}$-isotrivial $V$. Because we estimate $h(\psi(W))$ and not $h(W)$, it clearly suffices to assume that $\psi$ is the identity, so that $V$ is defined over ${\bf F}_q$. Take arbitrary $\pi$ in $V(\sqrt{G})$ not in $V({\bf F}_K)$. Then either $(i)$ or $(iii)$ of the Proposition holds.
\medskip
If $(i)$ holds, then $(b)$ looks good with $e=0$ (and $\psi$ the identity); at least $\pi$ lies in some $W(\sqrt{G})$ for a proper subvariety $W$ of $V$, defined over $K$, with
$$h(W) ~\leq~ 8n^54^nd\delta(n+r)R(\sqrt{G})^2. \eqno(9.1)$$
\medskip
What if $(iii)$ holds? Now any $a$ in ${\bf F}_q$ has a unique $p$th root $a^{{1 \over p}}$ in ${\bf F}_q$, which is also a conjugate of $a$ over ${\bf F}_p$. We get a new point $\pi'$ in $V'(\sqrt{G})$, also not in $V'({\bf F}_K)$, for a new variety $V'$ in ${\bf P}_n$ which is a conjugate of $V$. The new variety has the same dimension as $V$, and is also defined over ${\bf F}_q$. So we can repeat the process, and again we get either $(i)$ or $(iii)$ of the Proposition.
\medskip
If $(i)$ holds, then $\pi'$ lies in some $W'(\sqrt{G})$ again with $W'$ over $K$ and $h(W')$ bounded as in (9.1). So $\pi$ lies in $(W'(\sqrt{G}))^p$ as in (b) with $e=1$.
\medskip
Or if $(iii)$ holds, then we get a new point $\pi''$ in $V''(\sqrt{G})$ for a new conjugate $V''$ of $V$ in ${\bf P}_n$.
\medskip
And so on, in a manner similar to the looping in the $p$-automata of [D] section 4. Because $\pi$ was not in $V({\bf F}_K)$, this procedure must eventually stop at some proper subvariety $W^{(L)}$ over $K$ of $V^{(L)}$ (here the number $L$ of repetitions might depend on $\pi$). Now the original point $\pi$ lies in $(W^{(L)}(\sqrt{G}))^{p^L}$ with $h(W^{(L)})$ bounded as in (9.1).
\medskip
Because $\pi$ was arbitrary in $V(\sqrt{G})$ not in the finite set $V({\bf F}_K)$, the conclusion so far is 
$$V(\sqrt{G}) ~\subseteq~ \bigcup_{W \in {\cal W}}~\bigcup_{L=0}^\infty(W(\sqrt{G}))^{p^L} $$
for a collection ${\cal W}$ of proper subvarieties $W$ of conjugates of $V$ defined over $K$ and satisfying (9.1); here we may have to include single points $W$ with $h(W)=0$. To get equality we write $q=p^f$ and $L=fe+l$ for $e \geq 0$ and $0 \leq l \leq f-1$; this gives 
$$V(\sqrt{G}) ~\subseteq~ \bigcup_{\tilde W \in \tilde{\cal W}}~\bigcup_{e=0}^\infty(\tilde W(\sqrt{G}))^{q^e} $$
with a new collection $\tilde {\cal W}$ of proper subvarieties $\tilde W=W^{p^l}$ of conjugates of $V$ with
$$h(\tilde W)~=~p^lh(W) ~\leq~ 8n^54^n(q/p)d\delta(n+r)R(\sqrt{G})^2.$$
Finally by intersecting each $\tilde W$ with $V=V^q$ we can assume that each $\tilde W$ is a proper subvariety of $V$ itself in the above, without increasing the height further. Because $V$ is defined over ${\bf F}_q$, the $(\tilde W(\sqrt{G}))^{q^e}$ now lie in $(V(\sqrt{G}))^{q^e}=V(\sqrt{G})$, and so at last the two sides are equal. Now we have the desired $(b)$; of course the finiteness of the collection of $\tilde W$ follows from the Northcott Property already noted in section 7. This settles the case of transversal $\sqrt{G}$-isotrivial $V$ not contained in a proper coset.
\medskip
Henceforth (until further notice) we will assume that $V$ is not $\sqrt{G}$-isotrivial (and still transversal not contained in a proper coset). 
\medskip
Suppose first that $V$ is a hyperplane. Take arbitrary $\pi$ in $V(\sqrt{G})$. Then either $(i)$ or $(ii)$ of the Proposition holds. We regard this dichotomy as the starting stage $l=1$.
\medskip
If $(i)$ holds, then as before $(a)$ of the Main Estimate looks good; at least $\pi$ lies in some $W(\sqrt{G})$ for a proper subvariety $W$ of $V$, defined over $K$, with
$$h(W) \leq Ch(V)^{2n} \eqno(9.2)$$
for 
$$C~=~8n^54^nd\delta(n+r)R(\sqrt{G})^2.\eqno(9.3)$$
\medskip
What if $(ii)$ holds? We get a new point $\pi'$ in $V'(\sqrt{G})$ for a new variety $V'$ in ${\bf P}_n$ with 
$$\pi=\psi(\pi'^p),~~V=\psi(V'^p).\eqno(9.4)$$
Here $\psi$ is a $\sqrt{G}$-automorphism with
$$h(\psi) \leq pB \eqno(9.5)$$
for 
$$B=n\delta(n+r)R(\sqrt{G})^2. \eqno(9.6)$$
\medskip
This $V'$ is also a hyperplane, and also not $\sqrt{G}$-isotrivial. So we can repeat the process, and again we get either $(i)$ or $(ii)$ of the Proposition. This dichotomy is the next stage $l=2$.
\medskip
If $(i)$ holds, then $\pi'$ lies in some $W'(\sqrt{G})$. So $\pi$ lies in $W(\sqrt{G})$ for $W=\psi(W'^p)$,
almost as good as above, except that $h(W)$ could be larger than before. We take care of this later.
\medskip
Or if $(ii)$ holds, then we get a new point $\pi''$ in $V''(\sqrt{G})$ for a new variety $V''$ in ${\bf P}_n$.
\medskip
And so on. At stage $l$ we get either $\pi^{(l-1)}$ in a proper subvariety $W^{(l-1)}$ of $V^{(l-1)}$ with
$$h(W^{(l-1)}) ~\leq~ Ch(V^{(l-1)})^{2n} \eqno(9.7)$$
as in (9.2) and (9.3), or a new point $\pi^{(l)}$ in $V^{(l)}(\sqrt{G})$ for a new variety $V^{(l)}$ with 
$$\pi^{(l-1)}=\psi^{(l-1)}((\pi^{(l)})^p),~~V^{(l-1)}=\psi^{(l-1)}((V^{(l)})^p).\eqno(9.8)$$
as in (9.4), for
$$h(\psi^{(l-1)})  \leq pB. \eqno(9.9)$$
as in (9.5) and (9.6).
\medskip
We claim that this procedure must eventually stop because $V$ is not $\sqrt{G}$-isotrivial, and after a certain number $L$ of repetitions which this time is independent of $\pi$. Actually let us define the integer $L_0 \geq 0$ by
$$p^{L_0} ~\leq~ 2h(V)R(\sqrt{G})~<~p^{L_0+1}.\eqno(9.10)$$
From (9.8) we obtain $V=\psi_l((V^{(l)})^{p^l})$ with the $\sqrt{G}$-automorphism 
$$\psi_l~=~\psi\psi'^p\cdots(\psi^{(l-1)})^{p^{l-1}}.\eqno(9.11)$$
Writing the hyperplane $V$ in the affine form (8.1), we know that some coefficient $x=a_j \neq 0$ does not lie in $\sqrt{G}$, and $x=gy^{p^l}$ for some $g$ in $\sqrt{G}$ and some $y$ in $K$. We can now apply Lemma 4.5, because $\sqrt{G_k}$ there is just $\sqrt{G}$. We conclude that
$$p^l ~\leq~ 2h(x)R(\sqrt{G}) ~\leq~ 2h(V)R(\sqrt{G}).$$
In view of (9.10) this means that $(ii)$ cannot hold for $l=L_0+1$. Thus there is some $L$ with $0 \leq L \leq L_0$ such that $(ii)$ holds at stages $l=1,\ldots,L$ (at least if $L \geq 1$), and then $(i)$ holds at stage $l=L+1$. We conclude that $\pi^{(L)}$ lies in $W^{(L)}$, and from (9.7)
$$h(W^{(L)}) ~\leq~ Ch(V^{(L)})^{2n}. \eqno(9.12)$$
Thus $\pi=\psi_L((\pi^{(L)})^{p^L})$ lies in $W=\psi_L((W^{(L)})^{p^L})$. By (7.1) and (9.11) we get
$$h(W) ~\leq~ p^Lh(W^{(L)})+nh(\psi_L) ~\leq~ p^Lh(W^{(L)})+n\left(h(\psi)+ph(\psi')+\cdots+p^{L-1}h(\psi^{(L-1)}\right),$$
which using (9.9) and (9.12) yields
$$h(W) ~\leq~ Cp^Lh(V^{(L)})^{2n}+2np^LB ~\leq~ C(p^Lh(V^{(L)})^{2n}+2np^LB.\eqno(9.13)$$ 
To estimate $h(V^{(L)})$ we use (7.1), (9.8) and (9.9) to get
$$ph(V^{(l)})~=~h((\psi^{(l-1)})^{-1}V^{(l-1)}) ~\leq~ h(V^{(l-1)})+n^2h(\psi^{(l-1)})~\leq~ h(V^{(l-1)})+n^2pB.$$
If $L \geq 1$ we multiply this by $p^{l-1}$ and sum from $l=1$ to $l=L$, getting $p^Lh(V^{(L)}) \leq h(V)+2n^2p^LB$ (which holds also if $L=0$). Inserting this into (9.13) we get
$$h(W) ~\leq~ C\left(h(V)+2n^2p^LB\right)^{2n}+2np^LB~\leq~ 2C\left(h(V)+2n^2p^LB\right)^{2n},$$
and then using (9.6) and (9.10) with $L \leq L_0$ we find
$$h(W)~\leq~ 2Ch(V)^{2n}\left(1+4n^3\delta(n+r)R(\sqrt{G})^3\right)^{2n} ~\leq~ 2Ch(V)^{2n}\left(5n^3\delta(n+r)R(\sqrt{G})^3\right)^{2n}$$
From (9.3) we get finally 
$$h(W) ~\leq~C'h(V)^{2n}R(\sqrt{G})^{6n+2}\eqno(9.14)$$
with
$$C'~=~16n^54^nd\delta(n+r)\left(5n^3\delta(n+r)\right)^{2n} ~\leq~ 2n^2d\left(10n^3\delta(n+r)\right)^{2n+1}.$$
Because $\pi$ was arbitrary, the conclusion so far is 
$$V(\sqrt{G}) ~\subseteq~ \bigcup_{W \in {\cal W}}W(\sqrt{G})$$ 
for a finite collection ${\cal W}$ of proper subvarieties $W$ of $V$ satisfying (9.14). But then the two sides are of course equal.  This settles the Main Estimate for transversal hyperplanes $V$ that are not $\sqrt{G}$-isotrivial and not contained in a proper coset. 
\medskip
Next suppose that $V$, still not $\sqrt{G}$-isotrivial (and still transversal not contained in a proper coset), has dimension $m-1$ for some $m<n$. So after a permutation of variables it can be defined by equations (6.1). Each of these equations defines a hyperplane $V_i$, so that $V=V_m \cap \cdots \cap V_n$. 
\medskip
We claim that we can assume that all non-zero $a(i,j)$ lie in $\sqrt{G}$. Otherwise for example $V_m$ is transversal and not $\sqrt{G}$-isotrivial in the projective space with coordinates $X_j$ corresponding to $j=m$ and the $j$ with $a(m,j) \neq 0$. Since no $X_m-aX_j~(m \neq j,a \neq 0)$ vanishes on $V$, this projective space has dimension at least 2. So then we could apply the hyperplane result (9.14) to deduce that all solutions lie in a finite union of proper subspaces $W_m$ of this $V_m$ with
$$h(W_m) ~\leq~ C'h(V_m)^{2n}R(\sqrt{G})^{6n+2}.$$
But as in the affine situation just after (8.13), it can be seen that $W_m$ intersects the remaining intersection $U_m=\bigcap_{i \neq m}V_{i}$ in a proper subspace of $V=V_m \cap U_m$. For example the triangular nature of (6.1) makes it clear that $X_{m+1},\ldots,X_n$ are determined by $X_0,\ldots,X_{m-1}$ on $U_m$, and then that $X_m$ is determined by $X_0,\ldots,X_{m-1}$ on $W_m$ in $V_m$; but also some non-zero linear form in $X_0,\ldots,X_{m-1}$ must vanish on $W_m$. Therefore $W=W_m\cap U_m$ has dimension strictly less than $m-1$. So we are indeed in a proper subspace as required by (a) of the Main Estimate. Further $W=W_m \cap V$ and so $h(W) \leq h(W_m)+h(V)$ by Lemma 7.1; moreover $h(V_m) \leq h(V)$ because the $a(m,j)$ are themselves among the Grassmannian coordinates of $V$. We end up with (9.14) with say an extra factor 2.
\medskip
So indeed from now on we can assume that all non-zero $a(i,j)$ in (6.1) lie in $\sqrt{G}$. This means that we are set up to apply Lemma 6.1. We will see that the effect is to pass to a proper subvariety of at least one of $V_m,\ldots,V_n$ despite their being separately isotrivial. As $V$ is not $\sqrt{G}$-isotrivial by assumption, we find some quotient (6.2), say $Q$, not lying in ${\bf F}_K$. Let $\pi=(\xi_0,\ldots,\xi_n)$ be any point of $V(\sqrt{G})$. For a typical factor ${a(i,j) \over a(i,j')}$ in $Q$ we apply part $(b)$ of the Main Estimate in lower dimensions to $V_i$, with $\psi_i$ determined by 1 and the  non-zero $a(i,j)$. So here $q=p$. We find finitely many proper subspaces $W_i$ of $V_i$ such that $\psi_i(V_i(\sqrt{G}))$ lies in the union of the $\bigcup_{e=0}^\infty(\psi_i(W_i)(\sqrt{G}))^{p^e}$, with
$$h(\psi_i(W_i)) ~\leq~8n^54^nd\delta(n+r)R(\sqrt{G})^2 \eqno(9.15)$$ 
(now independent of $p$). In particular, writing $\pi_i$ for the projection of $\pi$ to the lower dimensional space, we have equations 
$$\psi_i(\pi_i)=\sigma_i^{q_i} \eqno(9.16)$$ 
for $\sigma_i$ in some $\psi_i(W_i)$ and some power $q_i$ of $p$. Thus
${a(i,j)\xi_{j} \over a(i,j')\xi_{j'}}=\eta^{q_i}$ for certain $\eta=\eta(i,j,j')$ in $K^*$. Multiplying all these over the factors in (6.2) we find $Q=\eta_1^{q_1} \cdots \eta_k^{q_k}$ for certain $\eta_1,\ldots,\eta_k$ in $K^*$. Because the fixed $Q$ is not in ${\bf F}_K$, this forces $q=\min\{q_1,\ldots,q_k\}$ to be bounded above by some quantity depending only on $V$. In fact $h(Q) \geq q$, but on the other hand from (6.2) we see that $h(Q) \leq (n+1)h(V)$. Thus 
$$q \leq (n+1)h(V). \eqno(9.17)$$
\medskip
Say this minimum is $q=q_i$. Now (9.16) says that $\pi_i$ and so $\pi$ lies in the variety $U=\psi_i^{-1}(\psi_i(W_i))^{q}$ of dimension strictly less than the dimension of $V_i$. This intersects $V_i$ in a proper subvariety $W_i'$ of $V_i$. Once more this $W_i'$ intersects the remaining intersection $\bigcap_{i' \neq i}V_{i'}$ in a proper subvariety $W$ of $V$. As for heights, we have $W=W_i' \cap V$ so $h(W) \leq h(W_i')+h(V)$. Also $h(W_i') \leq h(U)+h(V_i) \leq h(U)+h(V)$, and also
$$h(U) ~\leq~ qh(\psi_i(W_i))+nh(\psi_i^{-1}) ~\leq~ qh(\psi_i(W_i))+n^2h(V_i)$$
because of the definition of $\psi_i$. Putting these together and using (9.15),(9.17) we conclude that
$$h(W) ~\leq~ 8n^5(n^2+n+3)4^nd\delta(n+r)h(V)R(\sqrt{G})^2.$$
This is much smaller than (9.14), and so we have completed the proof of the Main Estimate when $V$ is transversal and not contained in a proper coset. In case (a) we have reached so far the bound $h(W) \leq Ah(V)^{2n}R^{6n+2}$ with $R=R(\sqrt{G})$ and $A=4n^2d(10n^3\delta(n+r))^{2n+1}$ due to the extra factor 2 encountered after establishing (9.14).
\medskip
To treat the more general situation when $V$ is transversal and not itself a coset, we use induction on $n \geq 2$, and we will obtain in case (a) the slightly weaker result $h(W) \leq Ah(V)^{2n}R^{6n+2}+nh(V)$. This leads at once to the bound given in the Main Estimate.
\medskip
If $n=2$ then there is a single equation $a_0X_0+a_1X_1+a_2X_2=0$, and transversality implies all $a_i \neq 0$. Thus no $X_i-aX_j~(i \neq j,a \neq 0)$ vanishes on $V$, and we are done. Thus we can suppose that $n \geq 3$.
\medskip
After permuting the variables, we can suppose that $X_n-aX_{n-1}~(a \neq 0)$ vanishes on $V$. In the remaining equations for $V$ we may eliminate $X_n$ to obtain a linear variety $\tilde V$ in ${\bf P}_{n-1}$. This $\tilde V$ cannot be a coset otherwise $V$ would be. Also $\tilde V$ certainly involves the variables $X_0,\ldots,X_{n-2}$ and so is transversal in ${\bf P}_{\tilde n}$ for $\tilde n=n-2$ or $\tilde n=n-1$. Here $\tilde n \geq 2$ unless $n=3$; but in that case if $\tilde V$ is not transversal in ${\bf P}_2$ then $V$ would be defined by equations $X_3=aX_2$ and $b_0X_0+b_1X_1=0$ so would be a coset. Thus we can assume that $\tilde V$ is transversal in ${\bf P}_{\tilde n}$ with $\tilde n \geq 2$.
\medskip
Suppose first that $V$ is not $\sqrt{G}$-isotrivial as in (a). Then $\tilde V$ cannot be $\sqrt{G}$-isotrivial otherwise we could transform $X_n$ to make $V$ isotrivial. Thus by induction the Main Estimate holds for $\tilde V$. It is now relatively straightforward to deduce the Main Estimate for $V$. Thus by case (a) for $\tilde V$ we get 
$$\tilde V(\sqrt{G})~=~\bigcup_{\tilde W \in \tilde{\cal W}}\tilde W(\sqrt{G})\eqno(9.18)$$
for a finite set $\tilde{\cal W}$ of proper linear subvarieties $\tilde W$ of $\tilde V$, also defined over $K$ and with  $h(\tilde W) \leq Ah(\tilde V)^{2n}R^{6n+2}+(n-1)h(\tilde V)$. Now we will check that (a) for $V$ follows with $W$ defined by the equations of $\tilde W$ together with $X_n=aX_{n-1}$. First the upper bound of Lemma 7.1 gives 
$$h(W) ~\leq~ h(\tilde W)+h(a) ~\leq~ Ah(\tilde V)^{2n}R^{6n+2}+(n-1)h(\tilde V)+h(a). \eqno(9.19)$$
We can suppose $X_{n-1} \neq 0$ on $\tilde V$, else (9.18) would be empty; and so the lower bound of Lemma 7.1 gives $h(V) \geq \max\{h(\tilde V),h(a)\}$. Therefore (9.19) implies
$$h(W) ~\leq~ Ah(V)^{2n}R^{6n+2}+nh(V)$$
as required.
\medskip
And in case (b) for $\sqrt{G}$-isotrivial $V$ (assuming as above that $\psi$ is the identity) we see that $\tilde V$ is $\sqrt{G}$-isotrivial and $a$ lies in ${\bf F}_q$. We get (b) for $V$ from (b) for $\tilde V$ using the analogue $\tilde V(\sqrt{G})=\bigcup_{\tilde W \in \tilde{\cal W}}~\bigcup_{e=0}^\infty(\tilde W(\sqrt{G}))^{q^e}$ of (9.18) with as above $W$ defined by the equations of $\tilde W$ together with $X_n=aX_{n-1}$; now $h(W) \leq h(\tilde W)$.
\medskip
What if $V$ is not transversal (and of course still not a coset)? Then it is transversal (and still not a coset) in some projective subspace of dimension $n' \leq n-1$. Here $n' \geq 2$ otherwise it would be a coset. The above cases (a) and (b) in dimension $n'$ now lead immediately to the same cases in ${\bf P}_n$; we have merely ignored $n-n'$ projective variables that were never in the equations anyway. 
\medskip
This finally finishes the proof of the Main Estimate.
\bigskip 
In view of the fact that the estimate in case (a) is independent of the characteristic $p$, it may seem a nuisance that the estimate in case (b) depends on $p$. But actually this is unavoidable, and there are even examples to show that the full $q/p$ is needed. To see this, take any power $q>1$ of $p$, and define $K={\bf F}_q(t)$ with $G=\sqrt{G}$ generated by $t,1-t$ and a generator $\zeta$ of ${\bf F}_q^*$. Here we have $r=2,R(\sqrt{G})=\sqrt{3}$ and, with the obvious transcendence basis, $d=1$. The affine equations
$$x+y=1,~~x+\zeta z=1$$
give rise to a $\sqrt{G}$-isotrivial line $V$ (with $h(V)=0$ and $\psi$ the identity), and an upper bound $B$ in (b) would mean that all solutions over $\sqrt{G}$ are given by $w,w^q,w^{q^2},\ldots$ for some $w$ with $h(w) \leq B$. Thus every solution $\pi$ would have either $h(\pi) \leq B$ or $h(\pi) \geq q$. But 
$$\pi~=~(x,y,z)~=~\left((1-t)^{q/p},t^{q/p},{t^{q/p} \over \zeta}\right)$$
is a solution with $h(\pi)=q/p$. It follows that $B \geq q/p$.
\bigskip
\bigskip
\noindent
{\bf 10. Isotrivial $W$.} We show here how to ensure that all the subvarieties $W$ in the Main Estimate can be made $\sqrt{G}$-isotrivial, at the expense of enlarging the exponents in the upper bounds for their heights. To simplify the various expressions we abbreviate the factors in case (a) of the Main Estimate by
$$\Delta=\Delta(n,r,d)=8n^2d(10n^3\delta(n+r))^{2n+1}\geq 1,~~h=h(V),~~R=R(\sqrt{G}), \eqno(10.1)$$
and that in case (b) of the Main Estimate by
$$\Psi=\Psi(n,r,d,p,q)=8n^54^n(q/p)d\delta(n+r) \geq 1.\eqno(10.2)$$
We also define some exponents
$$\rho(m)=\rho_n(m)={(2n)^m-1 \over 2n-1},~~~\eta(m)=\eta_n(m)=(2n)^m~~~~(m=1,2,\ldots)$$
\bigskip
\noindent
{\bf Main Estimate for isotrivial W}. {\it Let $V$ be a linear subvariety of ${\bf P}_n$ defined over $K$ but not a coset, with dimension $m-1 \geq 1$. 

\noindent
(a) If $V$ is not $\sqrt{G}$-isotrivial, then
$$V(\sqrt{G})~=~\bigcup_{W \in {\cal W}}W(\sqrt{G})$$
for a finite set $\cal W$ of proper linear $\sqrt{G}$-isotrivial subvarieties $W$ of $V$, also defined over $K$ and with
$$h(W) ~\leq~ (\Delta R^{6n+2})^{\rho(m)}h^{\eta(m)} \eqno(10.3)$$

\noindent
(b) If $V$ is $\sqrt{G}$-isotrivial and $\psi(V)$ is defined over ${\bf F}_q$, then
$$V(\sqrt{G})~=~\psi^{-1}\left(\bigcup_{W \in {\cal W}}~\bigcup_{e=0}^\infty(\psi(W)(\sqrt{G}))^{q^e}\right)$$
for a finite set $\cal W$ of proper linear $\sqrt{G}$-isotrivial subvarieties $W$ of $V$, also defined over $K$ and with
$$h(\psi(W)) ~\leq~ (\Delta R^{6n+2})^{\rho(m-1)}(\Psi R^2)^{\eta(m-1)}.$$}
\bigskip
\noindent
{\it Proof.} We start with case (a), and now we can write the bound as
$$h(W) ~\leq~ \Delta h^{2n}R^{6n+2} \eqno(10.4)$$
with $W$ not necessarily $\sqrt{G}$-isotrivial. We show by induction on the dimension $m-1 \geq 1$ of $V$ that the increased bound
$$h(\tilde W) ~\leq~ (\Delta R^{6n+2})^{\rho(m)}h^{\eta(m)}\eqno(10.5)$$
as in (10.3) holds where now all the $\tilde W$ are $\sqrt{G}$-isotrivial.
\medskip
When $m=2$ then the $W$ are points and so automatically $\sqrt{G}$-isotrivial as long as $W(\sqrt{G})$ is non-empty.
\medskip
When $m \geq 3$ we are fine unless some $W$ is not $\sqrt{G}$-isotrivial. We observe that such a $W$ cannot be a coset $T$. For the latter is defined by finitely many $X_i=a_{ij}X_j~(a_{ij} \neq 0)$, and if $T(\sqrt{G})$ is non-empty then clearly each $a_{ij}$ lies in $\sqrt{G}$. But now it is easy to see that $T$ is $\sqrt{G}$-isotrivial after all. For example we can rewrite the equations as $a_iX_i=a_jX_j$ with $a_i,a_j$ in $\sqrt{G}$. Then we can set up an equivalence relation on $\{0,1,\ldots,n\}$ characterized by the equivalence of such $i,j$. And now we need change only the variables in the equivalence classes of cardinality at least 2 in order to trivialize $T$. 
\medskip
So by induction each of these $W$ satisfies
$$W(\sqrt{G})~=~\bigcup_{\tilde W \in \tilde {\cal W}}\tilde W(\sqrt{G})$$
with $\sqrt{G}$-isotrivial $\tilde W$ such that
$$h(\tilde W) ~\leq~ (\Delta R^{6n+2})^{\rho(m-1)}h(W)^{\eta(m-1)}.$$
Therefore all we have to do is substitute (10.4) into this. We find the upper bound (10.5) because
$$\rho(m-1)+\eta(m-1)=\rho(m),~~~2n\eta(m-1)=\eta(m).$$
\medskip
For case (b) we write the bound as
$$h(\psi(W)) \leq \Psi R^2 \eqno(10.6)$$
with $W$ not necessarily $\sqrt{G}$-isotrivial. If some $W$ is not $\sqrt{G}$-isotrivial, then neither is $\psi(W)$, and we can write 
$$\psi(W)(\sqrt{G})~=~\bigcup_{W^* \in {\cal W}^*} W^*(\sqrt{G})$$
with $\sqrt{G}$-isotrivial $W^*$ such that
$$h(W^*) ~\leq~ (\Delta R^{6n+2})^{\rho(m-1)}h(\psi(W))^{\eta(m-1)}. \eqno(10.7)$$
Now we can see (without induction) that the bound
$$h(\psi(\tilde W)) ~\leq~ (\Delta R^{6n+2})^{\rho(m-1)}(\Psi R^2)^{\eta(m-1)} \eqno(10.8)$$
holds, where now all the $\tilde W=\psi^{-1}(W^*)$ are $\sqrt{G}$-isotrivial. In fact just as above, all we have to do is substitute (10.6) into (10.7), and we find at once (10.8). This completes the proof.
\bigskip
\bigskip
\noindent
{\bf 11. Points over $G$.} We show here how to replace $V(\sqrt{G})$ and $W(\sqrt{G})$ in the Main Estimate by $V(G)$ and $W(G)$ at the expense of worsening the dependence on the regulator.  However we no longer insist that the $W$ are isotrivial. If needed, this could be secured just by repeating the arguments of the previous section. We retain the notations (10.1),(10.2) from that section. Of course $n \geq 2$, and we continue with our assumption that $K$ is finitely generated over ${\bf F}_p$, with ${\bf F}_K=\overline{{\bf F}_p} \cap K$; further $G$ is finitely generated of rank $r \geq 1$ modulo ${\bf F}_K^*$.
\bigskip
\noindent
{\bf Main Estimate for points over $G$}. {\it There is a positive integer $f=f_K(G)\leq [\sqrt{G}:G]$, depending only on $K$ and $G$, with the following property. Let $V$ be a positive-dimensional linear subvariety of ${\bf P}_n$ defined over $K$ but not a coset. 

\noindent
(a) If $V$ is not $\sqrt{G}$-isotrivial, then 
$$V(G)~=~\bigcup_{W \in {\cal W}}W(G)$$
for a finite set $\cal W$ of proper linear subvarieties $W$ of $V$, also defined over $K$ and with
$$h(W) ~\leq~ \Delta h^{2n}R(\sqrt{G})^{6n+2}.$$

\noindent
(b) If $V$ is $\sqrt{G}$-isotrivial and $\psi(V)$ is defined over ${\bf F}_q$, then either

(ba) we have
$$V(G)~=~\bigcup_{W \in {\cal W}}W(G)$$
for a finite set $\cal W$ of proper linear subvarieties $W$ of $V$, also defined over $K$ and with
$$h(\psi(W)) ~\leq~ |{\bf F}_K|\Psi R(G)^2$$
or

(bb) we have
$$V(G)~=~\psi^{-1}\left(\bigcup_{W \in {\cal W}}~\bigcup_{e=0}^\infty(\psi(W)(G))^{q^{fe}}\right)\eqno(11.1)$$
for a finite set $\cal W$ of proper linear subvarieties $W$ of $V$, also defined over $K$ and with
$$h(\psi(W)) ~\leq~ q^{f} |{\bf F}_K|\Psi R(G)^2.\eqno(11.2)$$}
\bigskip
We need first a simple remark about congruences. Here $\phi$ is the Euler function.
\bigskip
\noindent
{\bf Lemma 11.1.} {\it For a given power $Q > 1$ of a prime $P$ consider a finite collection of congruence equations
$$LQ^e \equiv M {\rm ~mod} ~N \eqno(11.3)$$
with $N$ taken from a finite set ${\cal N}$ of positive integers and $L,M$ taken from ${\bf Z}$. Suppose that the set of solutions $e \geq 0$ is non-empty. Then if there is some $M \neq 0$ with ord$_PM <~$ord$_PN$ this set is

(a) finite with $Q^e \leq \max_{N \in {\cal N}}N$,

\noindent
and otherwise

(b) a finite union of arithmetic progressions $e=e_0,e_0+f,e_0+2f,\ldots$ with $f=\prod_{N \in {\cal N}}\phi(N)$ and $Q^{e_0} < Q^{f}\max_{N \in {\cal N}}N$.}
\bigskip
\noindent
{\it Proof.} Suppose first that there is some $M \neq 0$ with ord$_PM <~$ord$_PN$. Then the corresponding $L \neq 0$, and we get
$$e~ {\rm ord}_PQ ~\leq~ {\rm ord}_PLQ^e ~=~ {\rm ord}_PM ~<~ {\rm ord}_PN$$
giving case (a).
\medskip
Thus we can assume that ord$_PM \geq ~$ord$_PN$ whenever $M \neq 0$. We proceed to verify case (b). Now the congruences (11.3) can be split into congruences modulo powers of $P$ and congruences modulo powers $\tilde P^m$ of other primes $\tilde P \neq P$. 
\medskip
The former congruences, if any, will be satisfied as soon as $e$ is sufficiently large. Indeed they amount to $LQ^e \equiv 0~{\rm mod} ~P^{{\rm ord}_PN}$ and so conditions $e \geq \lambda$ for various real $\lambda \leq {{\rm ord}_PN \over {\rm ord}_PQ}$; that is, $Q^\lambda \leq P^{{\rm ord}_PN} \leq N$. Thus together they give a single condition $e \geq \Lambda$ for some real $\Lambda$ with $Q^\Lambda \leq \max_{N \in {\cal N}}N$. 
\medskip
We note that whether $e$ satisfies the other congruences depends only on its congruence class modulo $f$. For if $\tilde P^m$ divides some $N$ then $\phi(\tilde P^m)$ divides $\phi(N)$ which divides $f$, and so $Q^f \equiv 1$ mod $\tilde P^m$.
\medskip
Thus the solutions $e$ satisfy $e \geq \Lambda$ and also must lie in a finite number of arithmetic progressions modulo $f$. If $e_0$ is the smallest member of one of these progressions with $e_0 \geq \Lambda$, then $e_0-f < \Lambda$ and this leads to case (b), thereby completing the proof.
\bigskip
We can now start on the proof of the Main Estimate for points over $G$.
\medskip
Suppose first that $V$ is not $\sqrt{G}$-isotrivial. Then (a) of the Main Estimate gives 
$$V(\sqrt{G})~=~\bigcup_{W \in {\cal W}}W(\sqrt{G})$$ 
for $W$ satisfying (10.4). Now we can descend to $G$ simply by intersecting with ${\bf P}_n(G)$.
\medskip
Next suppose that $V$ is $\sqrt{G}$-isotrivial and $\psi(V)$ is defined over ${\bf F}_q$. Using elementary divisors we can find generators $\gamma_1,\ldots,\gamma_r$ of $\sqrt{G}$ modulo constants and positive integers $d_1,\ldots,d_r$ such that $\gamma_1^{d_1},\ldots,\gamma_r^{d_r}$ generate $G$ modulo constants. The constants can be taken care of with an extra $\gamma_0$ generating $\sqrt{G} \cap {\bf F}_K$ and $\gamma_0^{d_0}$ generating $G \cap {\bf F}_K$; here $d_0$ divides the order of $\gamma_0$ as a root of unity. Thus
$$[\sqrt{G}:G]=d_0d_1\cdots d_r. \eqno(11.4)$$
We write
$$\psi(X_0,\ldots,X_n)=(\psi_0X_0,\ldots,\psi_nX_n)$$
with 
$$\psi_i~=~\gamma_0^{a_{0i}}\gamma_1^{a_{1i}}\cdots \gamma_r^{a_{ri}}~~~(i=0,\ldots,n)\eqno(11.5)$$
in $\sqrt{G}$. Now (b) of the Main Estimate gives 
$$V(\sqrt{G})~=~\psi^{-1}\left(\bigcup_{W \in {\cal W}}~\bigcup_{e=0}^\infty(\psi(W)(\sqrt{G}))^{q^e}\right)\eqno(11.6)$$
for $W$ satisfying (10.6). But we can no longer descend to $G$ simply by intersecting with ${\bf P}_n(G)$.
\medskip
Consider a point $\pi=(\pi_0,\ldots,\pi_n)$ of $V(G)$. By (11.6) there is a point $\sigma=(\sigma_0,\ldots,\sigma_n)$ in some $W(\sqrt{G})$ and some $e \geq 0$ such that $\pi=\psi^{-1}(\psi(\sigma))^{q^e}$. As in (11.5) we write
$$\sigma_i~=~\gamma_0^{b_{0i}}\gamma_1^{b_{1i}}\cdots \gamma_r^{b_{ri}}~~~(i=0,\ldots,n);\eqno(11.7)$$ 
however $\pi$ is over $G$ and so
$$\pi_i~=~\gamma_0^{c_{0i}d_0}\gamma_1^{c_{1i}d_1}\cdots \gamma_r^{c_{ri}d_r}~~~(i=0,\ldots,n).$$
Equating exponents we find a system of congruences
$$(a_{ji}+b_{ji})q^e \equiv a_{ji}~~~{\rm mod} ~d_j~~~(i=0,\ldots,n;~j=0,1,\ldots,r)\eqno(11.8)$$
depending only on $\sigma$. We can apply Lemma 11.1, and the argument splits into two according to the conclusion. As the $b_{ji}$ in (11.7) appear only in the coefficients $L$, the splitting is independent of $\sigma$.
\medskip
Suppose first that Lemma 11.1(a) holds. Then 
$$q^e ~\leq~ \max\{d_0,d_1,\ldots,d_r\} ~\leq~ d_0d_1\cdots d_r ~=~ [\sqrt{G}:G]\eqno(11.9)$$ 
by (11.4). Now $\pi$ lies in the finitely many $\tilde W=\psi^{-1}(\psi(W))^{q^e}$, which we can put together into a set $\tilde{\cal W}$, and then we have shown that
$$V(G) ~\subseteq~ \bigcup_{\tilde W \in \tilde {\cal W}}\tilde W(\sqrt{G}).$$
Now intersecting with ${\bf P}_n(G)$ gives the same inclusion but with $\tilde W(G)$ on the right-hand side. On the other hand 
$$\tilde W~=~\psi^{-1}(\psi(W))^{q^e} ~\subseteq~ \psi^{-1}(\psi(V))^{q^e}~=~\psi^{-1}(\psi(V))~=~V$$ 
because $\psi(V)$ is defined over ${\bf F}_q$. Thus we conclude
$$V(G) ~=~ \bigcup_{\tilde W \in \tilde {\cal W}}\tilde W(G)$$
as in (ba) of the Main Estimate for points over $G$. But now from (11.9) and (10.6) the heights satisfy
$$h(\psi(\tilde W))~=~q^eh(\psi(W))~\leq~ d_0d_1\cdots d_r\Psi R(\sqrt{G})^2.$$
Using Lemma 4.1 we see that $R(G)=d_1\cdots d_rR(\sqrt{G})$, and so we can absorb some terms into the regulator to get
$$h(\psi(\tilde W))~\leq~ d_0\Psi R(G)^2 ~\leq~ |{\bf F}_K|\Psi R(G)^2. \eqno(11.10)$$
This completes the proof of (ba).
\medskip
It remains only to suppose that Lemma 11.1(b) holds. Then we know that $e=e_0+f\tilde e$ with $\tilde e \geq 0$ and $e_0$ bounded as in (11.9) but with an extra $q^{f}$. In particular taking $\tilde e=0$ we get a solution of (11.8) and this means that $\tilde \sigma=\psi^{-1}(\psi(\sigma))^{q^{e_0}}$ is also defined over $G$. It lies in
$$\tilde W=\psi^{-1}(\psi(W))^{q^{e_0}} \eqno(11.11)$$
and so in $\tilde W(G)$. We also have
$$\psi(\pi)=(\psi(\sigma))^{q^e}=(\psi(\tilde \sigma))^{{\tilde q}^{\tilde e}}$$
for $\tilde q=q^f$. Thus we conclude
$$V(G)~\subseteq~ \psi^{-1}\left(\bigcup_{\tilde W \in \tilde {\cal W}}~\bigcup_{\tilde e=0}^\infty(\psi(\tilde W)(G))^{{\tilde q}^{\tilde e}}\right)\eqno(11.12)$$
for the finite set $\tilde {\cal W}$ of $\tilde W$ in (11.11). On the other hand 
$$\psi(\tilde W)^{{\tilde q}^{\tilde e}}~=~(\psi(W))^{q^{e_0}{\tilde q}^{\tilde e}} ~\subseteq~ (\psi(V))^{q^{e_0}{\tilde q}^{\tilde e}}~=~\psi(V)$$ 
again because $\psi(V)$ is defined over ${\bf F}_q$. Thus we conclude equality in (11.12).
\medskip
Finally we calculate that $h(\psi(\tilde W))=q^{e_0}h(\psi(W))$ is bounded above by
$$q^{f} \max\{d_0,d_1,\ldots,d_r\}\Psi R(\sqrt{G})^2~\leq~ q^{f} |{\bf F}_K|\Psi R(G)^2\eqno(11.13)$$ 
as in (11.10), and of course $f =\phi(d_0)\phi(d_1)\cdots \phi(d_r)$ depends only on $K$ and $G$ with
$$f~\leq~ d_0d_1\cdots d_r~=~[\sqrt{G}:G].$$
This completes the proof of (bb); and so the Main Estimate for points over $G$ is proved.
\bigskip
In (11.13) the term $q^{f}$ cannot be so easily absorbed into the regulator without introducing an exponential dependence on $R(G)$. Let us discuss some aspects of this.
\medskip
When $G=\sqrt{G}$ then $f=1$ in (bb) and we are more or less back to (b) of the Main Estimate. But in general we need the extra $f$ in (11.1). The following example shows that it sometimes must be almost as large as $[\sqrt{G}:G]$. 
\medskip
We go back to the equation $t^mx+y=1$ of (1.5) over $K={\bf F}_p(t)$, with $n=2$. It is to be solved in the group $G=G_l$ generated by $t^l$ and $1-t$, so that $r=2$. Here $\sqrt{G}$ is generated by $t$ and $1-t$ together with a generator $\zeta$ of ${\bf F}_p^*$. The equation defines a $\sqrt{G}$-isotrivial line $V$ with $\psi(x,y)=(t^mx,y)=(\tilde x,\tilde y)$, so that $\tilde V=\psi(V)$ is defined by $\tilde x+\tilde y=1$, with $q=p$. 
\medskip
Now Leitner [Le] has found all points on $\tilde V(\sqrt{G})$. If $p$ is odd there are $p-2$ constant points in ${\bf F}_p^2$ together with six infinite families 
$$(\tilde x,\tilde y)=(\tilde x_{0}^{p^{\tilde e}},\tilde y_{0}^{p^{\tilde e}})~~~~(\tilde e=0,1,\ldots),$$
where $(\tilde x_{0},\tilde y_{0})$ are given by
$$(t,1-t),~(1-t,t),~\left({1 \over t},-{1-t \over t}\right),~\left(-{1-t \over t},{1 \over t}\right),~\left({1 \over 1-t},-{t \over 1-t}\right),~\left(-{t \over 1-t},{1 \over 1-t}\right).$$
The $(x,y)=\psi^{-1}(\tilde x,\tilde y)=(t^{-m}\tilde x,\tilde y)$ are all the points on $V(\sqrt{G})$. Choosing $m$ not divisible by $l$, we see that none of the constant points give rise to points of $V(G)$. Similarly for the second family above. And the same is true of the last four families above, simply because of the minus signs. However the first family gives $(t^{-m}t^{p^{\tilde e}},(1-t)^{p^{\tilde e}})$, which is in $G^2$ if and only if 
$$p^{\tilde e} \equiv m~ {\rm mod}~ l.\eqno(11.14)$$
\medskip
Now Artin's Conjecture implies that given any prime $p$, there are infinitely many primes $l$ for which $p$ is a primitive root modulo $l$. And Heath-Brown's Corollary 2 of [He] (p.27) implies that this is true for at least one of $p=3,5,7$. We can choose $m$ with $1 \leq m <l$ with $p^{l-2} \equiv m$ mod $l$. Now (11.14) implies $\tilde e \equiv l-2$ mod $l-1$ so $\tilde e=l-2+(l-1)e ~(e=0,1,\ldots)$. Thus the surviving points on $V(G)$ are just the 
$$\pi~=~\psi^{-1}(\psi(W))^{p^{(l-1)e}} ~~~~(e=0,1,\ldots)\eqno(11.15)$$ 
with $W$ as the single point $(t^{-m}t^{p^{l-2}},(1-t)^{p^{l-2}})$. This makes it clear that $f \geq l-1$ in (11.1); almost as big as $[\sqrt{G}:G]=(p-1)l$ for fixed $p$.
\medskip
We could also see this from (11.2). For as $R(G)=l\sqrt{3}$, it implies that there would be a point $\pi$ on $V(G)$ with $h(\psi(\pi)) \leq cp^fl^2$ for $c$ absolute. But the point (11.15) has $y=\tilde y=(1-t)^{p^{l-2}p^{(l-1)e}}$ so 
$$h(\psi(\pi)) ~\geq~ p^{l-2}p^{(l-1)e}~\geq~ p^{l-2}.\eqno(11.16)$$
Making $l \to \infty$, we deduce $f \geq l-c'\log l$, also almost as big as $[\sqrt{G}:G]=(p-1)l$.
\medskip
Less precisely, there can be no estimate
$$h(\psi(W)) ~\leq~ C(n,r,K)\left(h(V)+R(G)\right)^{\kappa}$$
replacing (11.2) which is polynomial in $h(V)$ and $R(G)$ for fixed $n,r,K$. For this would give a point with $h(\psi(\pi)) \leq c''(m+l)^{\kappa}\leq c'''l^{\kappa}$, contradicting (11.16). Similarly one sees that if the dependence on $h(V)$ is polynomial, then the dependence on $R(G)$ must be exponential. This explains the large solutions like (1.16), with $p=2,l=83,m=42$.
\bigskip
\noindent
{\bf 12. Proof of Descent Steps and Theorems.} In the Descent Steps the variety $V$ is certainly defined over a finitely generated transcendental extension $K$ of ${\bf F}_p$, and now we can choose any separable transcendence basis to obtain a height function. Now the Descent Step over $\sqrt{G}$ follows from the Main Estimate for isotrivial $W$. And the Descent Step over $G$ follows, at least without the assumption that the $W$ are $\sqrt{G}$-isotrivial, from the Main Estimate for points over $G$. This assumption can be removed by induction just as in section 10 (without bothering about estimates): any $W$ that is not $\sqrt{G}$-isotrivial can be replaced by a finite union of $\sqrt{G}$-isotrivial varieties. 
\medskip
To prove Theorem 1 we may assume that $V$ has positive dimension. We apply the Main Estimate for points over $G$ repeatedly, taking always $q=|{\bf F}_K|^{f_K(G)}$ for safety. With $V_0=V$, an arbitrary point $\pi$ of $V_0(G)$ is either a point of $W(G)$ for finitely many $W$ in $V_0$ with $\dim W \leq  \dim V-1$, or a point $\psi_1^{-1}\varphi^{e_1}\psi_1(\pi_1)$ for $\pi_1$ in $V_1(G)$ for finitely many $V_1$ in $V_0$ with $\dim V_1 \leq  \dim V-1$ and some $e_1 \geq 0$, with $\psi_1(V_0)$ defined over ${\bf F}_K$. Then we argue similarly with $\pi_1$; and so on. After at most $\dim V \leq n-1$ steps we descend to cosets $T=V_h$, and only finitely many $\psi_1,\ldots,\psi_h$ turn up on the way, leading to expressions as in (1.12) and thereby establishing Theorem 1.
\medskip
For later use we note that not just the varieties $T$ but also the whole unions $[\psi_1,\ldots,\psi_h]T$ lie in the variety $V$. Why is this? Well, a typical point of the union has the shape $\pi=(\psi_1^{-1}\varphi^{e_1} \psi_1) \cdots (\psi_h^{-1}\varphi^{e_h} \psi_h)(\tau)$ for some $e_1,\ldots,e_h$ and some $\tau$ in $T$. The descent for Theorem 1 provides linear varieties $V=V_0,V_1,\ldots,V_h=T$. Now clearly $\tau$ lies in $T$ inside $V_{h-1}$, so $\psi_h^{-1}\varphi^{e_h} \psi_h(\tau)$ lies in
$$\psi_h^{-1}\varphi^{e_h} \psi_h(V_{h-1})~=~\psi_h^{-1}\psi_h(V_{h-1})~=~V_{h-1}$$
inside $V_{h-2}$. In the same way $(\psi_{h-1}^{-1}\varphi^{e_{h-1}}\psi_{h-1})(\psi_h^{-1}\varphi^{e_h} \psi_h)(\tau)$ lies in $V_{h-2}$ inside $V_{h-3}$. Continuing backwards we see that $\pi=(\psi_1^{-1}\varphi^{e_1} \psi_1) \cdots (\psi_h^{-1}\varphi^{e_h} \psi_h)(\tau)$ lies in $V$.
\medskip
We leave it to the reader to check, by a straightforward induction argument like that in section 10 and also using Lemma 7.2, that for Theorem 1 one can take
$$\max\{h(\psi_1),\ldots,h(\psi_h),h(T)\}~\leq~(2q^2\Delta R(G)^{6n+2})^{\rho(m)}h(V)^{\eta(m)} \eqno(12.1)$$
in the notation of section 10. This indeed looks polynomial in $R(G)$ and $h(V)$; however, as we noted, an exponential dependence on $R(G)$ may be hiding in $q=|{\bf F}_K|^{f_K(G)}$.
\medskip
For the symmetrization argument in the proof of Theorem 2 we need a version of Lemma 8.1 (p. 209) of [D], partly removed from its recurrence context.
\bigskip
\noindent
{\bf Lemma 12.1.} {\it For $m \geq 1$ and $x_1,\ldots,x_m,y_1,\ldots,y_m$ in $K$ suppose that 
$$x_1y_1^{q^l}+\cdots+x_my_m^{q^l}~~=~~0 \eqno(12.2)$$ 
for all large $l$. Then this holds for all $l \geq 0$.}
\bigskip
\noindent
{\it Proof.} The proof will be by induction on $m$, the case $m=1$ being trivial. For the induction step we can clearly assume that $x_1,\ldots,x_m$ are non-zero. Now we note that (12.2) for any $m$ consecutive integers $l=g,g+1,\ldots,g+m-1$ implies the linear dependence of $y_1,\ldots,y_m$ over ${\bf F}_q$. For if we regard these as linear equations for $x_1,\ldots,x_m$, the underlying determinant is the $q^g$ power of that with entries $y_i^{q^{j-1}} (i,j=1,\ldots,m)$, and it is well-known that the latter, a so-called Moore determinant, is up to a constant the product of the $\beta_1y_1+ \cdots+\beta_my_m$ taken over all $(\beta_1,\ldots,\beta_m)$ in ${\bf P}_{m-1}({\bf F}_q)$ (see for example [Go] Corollary 1.3.7 p.8). Thus after permuting we can suppose that $y_m=\alpha_1y_1+\cdots+\alpha_{m-1}y_{m-1}$ for $\alpha_1,\ldots,\alpha_{m-1}$ in ${\bf F}_q$. Substituting into (12.2) gives
$$(x_1+\alpha_1x_m)y_1^{q^l}+\cdots+(x_{m-1}+\alpha_{m-1}x_m)y_{m-1}^{q^l}~=~0,$$
which therefore also holds for all large $l$. By the induction hypothesis we conclude that this holds for all $l \geq 0$, which leads back to (12.2) for all $l \geq 0$ and thus completes the proof.
\bigskip
To prove Theorem 2 consider a single $[\psi_1,\ldots,\psi_h]T(G)$ coming from Theorem 1. Fix $\tau_0$ in $T(G)$; then $T=\tau_0S$ for a linear subgroup $S$. 
\medskip
We argue first on the geometric level. According to (1.12) a typical point of $[\psi_1,\ldots,\psi_h]T$ has the shape
$$\psi_1^{q_1-1}\psi_2^{q_1q_2-q_1}\psi_3^{q_1q_2q_3-q_1q_2}\cdots \psi_h^{q_1\cdots q_h-q_1\cdots q_{h-1}}(\tau_0\sigma)^{q_1\cdots q_h}$$
with $q_i=q^{e_i}~(i=1,\ldots,h)$ and $\sigma$ in $S$; here we are regarding the $\psi_i~(i=1,\ldots,h)$ as multiplication by points instead of automorphisms. This expression can be written as
$$\pi_0\pi_1^{q_1}\pi_2^{q_1q_2}\cdots \pi_{h-1}^{q_1\cdots q_{h-1}}\pi_h^{q_1\cdots q_h}\sigma^{q_1\cdots q_h} \eqno(12.3)$$
with
$$\pi_0=\psi_1^{-1},~\pi_1=\psi_2^{-1}\psi_1,\ldots,~\pi_{h-1}=\psi_h^{-1}\psi_{h-1}, ~\pi_h=\psi_h\tau_0. \eqno(12.4)$$
Now when we write $q^{l_i}=q_1\cdots q_i ~(i=1,\ldots,h)$ we certainly get a point of $(\pi_0,\pi_1,\ldots,\pi_h)S$ according to (1.14); but at the moment we have asymmetry $l_1 \leq \cdots \leq l_h$. We eliminate the inequalities here as in [D] (p.212).
\medskip
Let us start with the last inequality. We can write (12.3) as $\xi\eta^{q^{l}}$ with $\xi$ and $\eta$ independent of $l=l_{h}$. We already remarked that $[\psi_1,\ldots,\psi_h]T$ lies in $V$, so (12.3) does. Thus for each linear form $\cal L$ defining $V$ we have ${\cal L}(\xi\eta^{q^{l}})=0$ for all $l_1,\ldots,l_{h-1},l$ with $0 \leq l_1 \leq \cdots \leq l_{h-1} \leq l$. Fixing $l_1,\ldots,l_{h-1}$, we see from Lemma 12.1 that this equation for all large $l$ implies the same equation for all $l \geq 0$. Thus the inequality $l_{h-1} \leq l_h$ has indeed been eliminated. Similar arguments work for the other conditions, as is clear from the arguments of [D] (p.212) after equation (22). For example, the next step fixes $l_1,\ldots, l_{h-2},l_h$ but not $l=l_{h-1}$.
\medskip
Looking back at (12.3), we have therefore proved that all the points
$$\pi_0\pi_1^{r_1}\pi_2^{r_2}\cdots \pi_{h-1}^{r_{h-1}}\pi_h^{r_h}\sigma^{r_h}\eqno(12.5)$$
lie in $V$, where the integers $r_i=q^{l_i}~(i=1,\ldots,h)$ now range independently over all positive integral powers of $q$. This is the required symmetrization at the geometric level. 
\medskip
It actually shows that the entire $(\pi_0,\pi_1,\ldots,\pi_h)S$ lies in $V$. For a typical point of the former has the shape
$$\pi_0\pi_1^{r_1}\pi_2^{r_2}\cdots \pi_{h-1}^{r_{h-1}}\pi_h^{r_h}\tilde \sigma\eqno(12.6)$$
for $\tilde\sigma$ in $S$. And there is $\sigma$ in $S$ with $\sigma^{r_h}=\tilde\sigma$. This could be interpreted as something about the divisibility of group varieties; but for us it is just a simple consequence of the fact that $S$ is defined by equations $X_i=X_j$. And now (12.6) and (12.5) are equal.
\medskip
At the arithmetic level we claim that $(\pi_0,\pi_1,\ldots,\pi_h)S(G)$ lies in $V(G)$. In fact every point
$$\pi~=~\pi_0\pi_1^{r_1}\pi_2^{r_2}\cdots \pi_{h-1}^{r_{h-1}}\pi_h^{r_h}\eqno(12.7)$$
with asymmetry $r_1 \leq \cdots \leq r_h$ has the shape (12.3) (with all coordinates of $\sigma$ equal to 1). It therefore lies in $[\psi_1,\ldots,\psi_h]T(G)$ which is in turn contained in $V(G)$. In particular $\pi$ lies in ${\bf P}_n(G)$. But why does it continue to lie in ${\bf P}_n(G)$ when the asymmetry is lifted?
\medskip
Well, we can take $r_1=\cdots=r_h=1$ in (12.7) to see that the product
$$\pi_0\pi_1\cdots \pi_h \eqno(12.8)$$
lies in ${\bf P}_n(G)$. Then taking $r_1=\cdots=r_{h-1}=1, r_h=q$ we can deduce that $\pi_h^{q-1}$ lies in ${\bf P}_n(G)$. And taking $r_1=\cdots=r_{h-2}=1, r_{h-1}=r_h=q$ we deduce that $\pi_{h-1}^{q-1}$ lies in ${\bf P}_n(G)$. And so on, until we see that all of
$$\pi_1^{q-1},\ldots, \pi_h^{q-1} \eqno (12.9)$$
lie in ${\bf P}_n(G)$ (this was already remarked in section 1).
\medskip
And now if $r_1,\ldots,r_h$ are arbitrary integral powers of $q$ in (12.7) we can write
$$\pi~=~(\pi_0\pi_1\cdots \pi_h)\pi_1^{r_1-1}\cdots \pi_h^{r_h-1}$$
to see from (12.8) and (12.9) that indeed $\pi$ lies in ${\bf P}_n(G)$.
\medskip
Now any point of $(\pi_0,\pi_1,\ldots,\pi_h)S(G)$ by (12.5) has the form $\pi\sigma^{r_h}$ with $\pi$ as above and $\sigma$ in $S(G)$. It follows that $(\pi_0,\pi_1,\ldots,\pi_h)S(G)$ lies in $V(G)$ as claimed. 
\medskip
On the other hand, taking all coordinates of $\sigma$ as 1 in (12.3) shows that $[\psi_1,\ldots,\psi_h]\{\tau_0\}$ lies in $(\pi_0,\pi_1,\ldots,\pi_h)S(G)$. As we could have fixed $\tau_0$ arbitrarily in $T(G)$, we see that $[\psi_1,\ldots,\psi_h]T(G)$ lies in $(\pi_0,\pi_1,\ldots,\pi_h)S(G)$.
\medskip
It follows that $V(G)$ is indeed the union of the $(\pi_0,\pi_1,\ldots,\pi_h)S(G)$, which completes the proof of Theorem 2. We note for later use the fact, already observed, that each $(\pi_0,\pi_1,\ldots,\pi_h)S$ is contained in $V$.
\medskip
Here too we leave it to the reader to check using (12.1) that for Theorem 2 one can take
$$\max\{h(\pi_0),h(\pi_1),\ldots,h(\pi_h)\}~\leq~(n+1)(2q^2\Delta R(G)^{6n+2})^{\rho(m)}h(V)^{\eta(m)}. \eqno(12.10)$$
This follows quickly from (12.4) and the easy fact that any $T(G)$ contains $\tau_0$ with $h(\tau_0) \leq h(T)$.
\medskip
To prove part (1) of Theorem 3 we start from Theorem 1 with $V=H$. We first claim that if some $\pi$ in $H(G)$ lies in some $[\psi_1,\ldots,\psi_h]T(G)$ with $T$ not a single point then some (1.2) fails for $\pi$. To see this, note that if $T$ is not a single point, then there is a partition of $\{0,1,\ldots,n\}$ into proper subsets $I,J, \ldots$ such that $T$ is defined by the proportionality of the homogeneous coordinates $X_i~(i \in I)$, $X_j~(j \in J)$, and so on. We may suppose that $I$ contains $0$ and that the equations corresponding to $I$ are $g_iX_0=g_0X_i$ for $i$ in $I$. Consider the point $\tau_I$ in ${\bf P}_n$ whose coordinates $X_i=g_i$ for $i$ in $I$ but with all other coordinates zero. It also lies in $T$.
\medskip
Now $\pi=(\psi_1^{-1}\varphi^{e_1} \psi_1) \cdots (\psi_h^{-1}\varphi^{e_h} \psi_h)(\tau)$ for some $e_1,\ldots,e_h$ and some $\tau$ in $T$. From our remark following the proof of Theorem 1, we see that $\pi_I=(\psi_1^{-1}\varphi^{e_1} \psi_1) \cdots (\psi_h^{-1}\varphi^{e_h} \psi_h)(\tau_I)$ lies in $H$. Now $\tau$ and $\tau_I$ have the same coordinates $X_i~(i \in I)$. It follows that $\pi$ and $\pi_I$ have the same coordinates $X_i~(i \in I)$. Since the other coordinates of $\pi_I$ are zero, this means that (1.2) fails for $\pi$ as claimed.
\medskip
Therefore $H^*(G)$ is contained in a finite union of sets $[\psi_1,\ldots,\psi_h]\{\tau\}$. And each of these lies in $H(G)$. This proves part (1) of Theorem 3. 
\medskip
Part (2) follows in a similar way with the help of the remark after the proof of Theorem 2, with $\pi=\pi_0(\varphi^{l_1} \pi_1) \cdots (\varphi^{l_h} \pi_h)\sigma$ and $\pi_I=\pi_0(\varphi^{l_1} \pi_1) \cdots (\varphi^{l_h} \pi_h)\sigma_I$ for $\sigma_I$ defined by $X_i=1$ for $i$ in $I$ but with all other coordinates zero. This shows that we can restrict to single points $S$, and the proof is finished as above. We have therefore proved all of Theorem 3.
\medskip
It is easy to deduce explicit estimates for Theorem 3 as for Theorems 1 and 2. One obtains at once (12.1) (with $T$ replaced by $\tau$) and (12.10).
\bigskip
\bigskip
\noindent
{\bf 13. Limitation results.} We show here that for each $n \geq 2$ the bounds $h \leq n-1$ in Theorems 1 and 2 cannot always be improved; and also that if $p >2$ the $\psi_1,\ldots,\psi_h$ in Theorem 1 and the $\pi_0,\pi_1,\ldots,\pi_h$ in Theorem 2 cannot always be chosen over $G$.
\medskip
We start with $h \leq n-1$. Because Theorem 1 directly implies Theorem 2 and then Theorem 3, it will suffice to prove the analogous statements for Theorem 3. Also we have seen that each $[\psi_1,\ldots,\psi_h]\{\tau\}$ in Theorem 3(1) is contained in some $(\pi_0,\pi_1,\ldots,\pi_h)$ in Theorem 3(2). So it is enough to treat Theorem 3(2).
\medskip
This we do with the affine hyperplane
$$x_1+x_2-x_3-\cdots-x_n~=~1\eqno(13.1)$$
already mentioned. 
\medskip
We need a simple observation. For a prime $p$ let $R=R_p$ be the set of points $(1,r_1,\ldots,r_{n-1})$ as the integers $r_1,\ldots,r_{n-1}$ run through all powers of $p$ satisfying the asymmetry conditions that $r_i$ divides $r_{i+1} ~(i=1,\ldots, n-2)$ and also the extra conditions
$$r_{n-1} ~\neq~ r_{n-2},r_{n-2}+r_{n-3},\ldots,r_{n-2}+r_{n-3}+\cdots+r_1. \eqno(13.2)$$
\bigskip
\noindent
{\bf Lemma 13.1.} {\it The set $R$ does not lie in a finite union of proper subgroups of ${\bf Z}^n$.}
\bigskip
\noindent
{\it Proof.} We can actually disregard (13.2) because their failure would just add more to the finite union of proper subgroups. Now the falsity of the lemma would lead to an equation
$${\cal F}(p^{e_1},\ldots,p^{e_{n-1}})~=~0 \eqno(13.3)$$
holding for all non-negative integers $e_1,\ldots,e_{n-1}$, where ${\cal F}(y_1,\ldots,y_{n-1})$ is a finite product of polynomials
$${\cal A}~=~a_0+a_1y_1+a_2y_1y_2+\cdots+a_{n-1}y_1y_2\cdots y_{n-1}$$
corresponding to the proper subgroups of ${\bf Z}^n$ perpendicular to $(a_0,\ldots,a_{n-1}) \neq 0$. It is clear that each ${\cal A} \neq 0$ and so ${\cal F} \neq 0$. On the other hand it is easy to see that the points in (13.3) are Zariski-dense in ${\bf R}^{n-1}$. This contradiction proves the lemma.
\bigskip
Take as usual $K={\bf F}_p(t)$ and $G$ generated by $t$ and $1-t$. We proceed to exhibit many points on $H^*(G)$ with $H$ defined by (13.1).
\medskip
For integral powers $q_1,\ldots,q_{n-1}$ of $p$ define
$$r_1=q_{n-1},~r_2=q_{n-1}q_{n-2},\ldots,~r_{n-1}=q_{n-1}\cdots q_1$$
and
$$~~~~~~d_1~=~r_{n-1}-r_{n-2}-\cdots-r_{2}-r_1,$$
$$d_2~=~r_{n-1}-r_{n-2}-\cdots-r_{2},$$
down to
$$~~~~d_{n-2}~=~r_{n-1}-r_{n-2}~~~~~~~~~~~~~~~~~~~~~$$
and
$$d_{n-1}~=~r_{n-1}.~~~~~~~~~~~~~~~~~~~~~~~~~$$
Then
$$x_1=t^{d_1},~x_2=1-t^{d_{n-1}},~x_3=t^{d_{n-2}}-t^{d_{n-1}},\ldots,x_n=t^{d_{1}}-t^{d_{2}}\eqno(13.4)$$
certainly satisfy (13.1), so the point $\xi=(x_1,\ldots,x_n)$ lies in $H$. It is in $H(G)$ because
$$x_2=1-t^{r_{n-1}}=(1-t)^{r_{n-1}},$$
$$x_3=t^{d_{n-2}}(1-t^{r_{n-2}})=t^{d_{n-2}}(1-t)^{r_{n-2}},$$
and so on.
\medskip
This also leads to a multiplicative representation
$$\xi~=~\xi_1^{r_1}\cdots \xi_{n-1}^{r_{n-1}} \eqno(13.5)$$
of the point in (13.4), where
$$\xi_1~=~({1 \over t},1,1,1,1,\ldots, 1,1, {1-t \over t}),$$
$$\xi_2~=~({1 \over t},1,1,1,1,\ldots, 1, {1-t \over t},{1 \over t})$$
$$\xi_3~=~({1 \over t},1,1,1,1,\ldots, {1-t \over t},{1 \over t},{1 \over t})$$
down to
$$\xi_{n-2}~=~({1 \over t},1,{1-t \over t},{1 \over t},{1 \over t},\ldots,{1 \over t},{1 \over t},{1 \over t}),$$
but
$$\xi_{n-1}~=~(t,1-t,t,t,t,\ldots,t,t,t).~~~~~~~~$$
We can quickly check that $\xi_1,\ldots, \xi_{n-1}$ are multiplicatively independent. Namely, a relation
$$\xi_1^{a_1}\cdots \xi_{n-1}^{a_{n-1}}~=~(1,1,1,1,1,\ldots,1,1,1)$$
would lead to $a_{n-1}=0$ on examining the second components, then $a_{n-2}=0$ from the third components, and so on down to $a_1=0$.
\medskip
The case $n=3$ with $q_1=q,q_2=r$ is of course (1.11) or (1.13).
\medskip
We can see that (13.4) lies in $H^*(G)$ provided $(1,r_1,\ldots,r_{n-1})$ lies in $R$. For the various exponents of $t$ clearly satisfy $d_{n-1}>d_{n-2}>\cdots >d_2>d_1$. There is one more exponent 0; but $d_{n-1} \neq 0$ and from the definition of $R$ we also have $d_{n-2} \neq 0,\ldots,d_1 \neq 0$. Thus the exponents $d_{n-1},\ldots,d_1,0$ in (13.4) are distinct, and it is easy to see that there can be no vanishing subsum of $x_1,x_2,-x_3,\ldots,-x_n$ (in fact each of $d_{n-2} = 0,\ldots,d_1 = 0$ does lead to a vanishing subsum). We already remarked that (1.13) is in $H^*$ as long as $r \neq s$, that is $q_1 \neq 1$, that is $r_2 \neq r_1$ as in (13.2).
\medskip
Now we can prove as promised that $H^*(G)$ does not lie in a finite union of sets
$$\Pi~=~(\pi_0,\pi_1,\ldots,\pi_h)_q~=~\bigcup_{l_1=0}^\infty \cdots \bigcup_{l_h=0}^\infty\pi_0\pi_1^{q^{l_1}}\cdots \pi_h^{q^{l_h}} \eqno(13.6)$$
for some $q$ and points $\pi_0,\pi_1\ldots,\pi_h$ with $h < n-1$. The idea is to note that each $\Pi$ lies in a coset of ${\bf G}_{\rm m}^n$ of dimension at most $h \leq n-2$; whereas the points (13.5) have rank $n-1$.
\medskip
Accordingly we assume that $H^*(G)$ does lie in such a finite union and we shall reach a contradiction.
\medskip
Now for each element of $R$ the corresponding (13.5) lies in $H^*(G)$ so in some $\Pi$. This provides a partition of $R$ into a finite union of subsets $R_\Pi$. By Lemma 13.1 we will be through if we can prove that each $R_\Pi$ lies in a proper subgroup of ${\bf Z}^n$.
\medskip
Suppose for some $\Pi$ we are lucky in the sense that the corresponding $\pi_0$ in (13.6) is multiplicatively independent of $\xi_1,\ldots,\xi_{n-1}$. The corresponding
$$\pi_0^{-1}\xi~=~\pi_0^{-1}\xi_1^{r_1}\cdots \xi_{n-1}^{r_{n-1}}$$
all lie in the group generated by $\pi_1,\ldots,\pi_h$, and so the multiplicative rank of the various $\pi_0^{-1}\xi$ is at most $h \leq n-2$. Since $\pi_0^{-1},\xi_1,\ldots, \xi_{n-1}$ are independent, it follows that the set $R_\Pi$ cannot contain $n$ (or even $n-1$) independent elements. So it must indeed lie in a proper subgroup of ${\bf Z}^n$.
\medskip
In fact we are not so likely to be that lucky, and it is more probable that there is a relation $\pi_0^a=\xi_1^{a_1}\cdots \xi_{n-1}^{a_{n-1}}$ with $a \neq 0$. Now the
$$\pi_0^{-a}\xi^a~=~\xi_1^{ar_1-a_1}\cdots \xi_{n-1}^{ar_{n-1}-a_{n-1}}$$
still lie in a group of rank at most $n-2$. Since $\xi_1,\ldots, \xi_{n-1}$ are independent, we deduce as above that the set of all $(ar_1-a_1,\ldots,ar_{n-1}-a_{n-1})$ lie in a proper subgroup of ${\bf Z}^{n-1}$. And this implies as above that $R_\Pi$ lies in a proper subgroup of ${\bf Z}^n$. 
\medskip
That finishes the proof of the first limitation result. We could also have argued with a symmetrized version of $R$; then the $\cal A$ in the proof of Lemma 13.1 could be taken more simply as $a_0+a_1y_1+a_2y_2+\cdots+a_{n-1}y_{n-1}$.
\medskip
We can use similar arguments to prove the second limitation result concerning non-definability over $G$. Because the $[\psi_1,\ldots,\psi_h]T(G)$ in Theorem 1 lead to $(\pi_0,\pi_1,\ldots,\pi_h)$ in Theorem 2 with (12.4) for $\tau_0$ in $T(G)$, it will again suffice to check the matter for Theorem 3(2).
\medskip
This we do with the affine line $H$ defined by $tx+y=1$ also over $K={\bf F}_p(t)$, now with $G$ generated by $t^{p-1}$ and $1-t$. It is the example treated at the end of section 11 with $m=1$ and $l=p-1$. We need another simple observation.
\bigskip
\noindent
{\bf Lemma 13.2.} {\it For an odd prime $p$ suppose that
$$q_1+q_2+q_3~=~\tilde q_1+\tilde q_2+\tilde q_3 \eqno(13.7)$$
for integral powers $q_1,q_2,q_3,\tilde q_1,\tilde q_2,\tilde q_3$ of $p$. Then $\tilde q_1,\tilde q_2,\tilde q_3$ are a permutation of $q_1,q_2,q_3$.}
\bigskip
\noindent
{\it Proof.} If $q_1,q_2,q_3$ are all different then the left-hand side of (13.7) has just three ones in its expansion to base $p$. So also the right-hand side; which means that $\tilde q_1,\tilde q_2,\tilde q_3$ are also all different. The result in this case is now clear (even for $p=2$). If say $q_1\neq q_2=q_3$ then we get a one and a two in the expansion because $p \neq 2$; so after a permutation $\tilde q_1\neq \tilde q_2=\tilde q_3$ too, and the result is still clear. Similarly if $q_1=q_2=q_3$ as long as $p \neq 3$. This last case can also be checked directly when $p=3$ and this proves the lemma; however the example $1+1+4=2+2+2$ shows that $p=2$ is not to be saved.
\bigskip
Now the analysis in section 11 before the primitive root business shows easily that the points of $H^*(G)=H(G)$ are given by
$$x=t^{r-1},~~y=(1-t)^r~~~~(r=1,p,p^2,\ldots). \eqno(13.8)$$
This is $(x,y)=\xi_0\xi_1^r$ for $\xi_0=(t^{-1},1)$ and $\xi_1=(t,1-t)$. Assume $p \neq 2$. If $H^*(G)$ were contained in a finite union of
$$\Pi~=~(\pi_0,\pi_1)_q~=~\bigcup_{l=0}^\infty \pi_0\pi_1^{q^{l}}$$
for some $q$ and some $\pi_0,\pi_1$ over $G$, then one of these $\Pi$ would certainly contain at least three different points (13.8). This gives equations
$$\xi_0\xi_1^r=\pi_0\pi_1^s,~~\xi_0\xi_1^{r'}=\pi_0\pi_1^{s'},~~\xi_0\xi_1^{r''}=\pi_0\pi_1^{s''} \eqno(13.9)$$
for powers $r<r'<r''$ of $p$ and powers $s,s',s''$ of $q$. Eliminating $\pi_0,\pi_1$ leads to
$$(\xi_0\xi_1^r)^{s'-s''}(\xi_0\xi_1^{r'})^{s''-s}(\xi_0\xi_1^{r''})^{s-s'}~=~1;$$
that is, $\xi_1^a=1$ for
$$a~=~r(s'-s'')+r'(s''-s)+r''(s-s').$$
So $a=0$; that is, 
$$rs'+r's''+r''s~=~rs''+r's+r''s'.$$
Lemma 13.2 shows in particular that $rs'$ is one of the terms on the right. But which one? Certainly $rs' \neq r''s'$. And $rs' \neq rs''$ else $s'=s''$ and (13.9) would imply $r'=r''$. It follows that $rs'=r's$. But now eliminating $\xi_1$ from the first two equations in (13.9) leads to $\xi_0^{r'-r}=\pi_0^{r'-r}$. Thus there would be $\alpha,\beta$ in $\overline{{\bf F}_p}$ with
$(\alpha t^{-1},\beta)=(\alpha,\beta)\xi_0=\pi_0;$
however this is impossible because $\alpha t^{-1}$ is not in $G$ if $p \neq 2$.
\bigskip
\bigskip
\bigskip
\noindent
\centerline{\bf References}
\bigskip
\noindent
[AB] B. Adamczewski and J.P. Bell, {\it On vanishing coefficients of algebraic power series over fields of positive characteristic}, Manuscript 2010.
\medskip
\noindent
[AV] D. Abramovich and F. Voloch, {\it Toward a proof of the Mordell-Lang conjecture in characteristic $p$}, International Math. Research Notices {\bf 5} (1992), 103-115.
\medskip
\noindent
[BG] E. Bombieri and W. Gubler,  {\it Heights in diophantine geometry}, New Mathematical Monographs {\bf 4}, Cambridge 2006. 
\medskip
\noindent
[BMZ] E. Bombieri, D. Masser and U. Zannier, {\it Intersecting a plane with algebraic subgroups of multiplicative groups}, Ann. Scuola Norm. Sup. Pisa Cl. Sci. (5) VII (2008), 51-80.
\medskip
\noindent
[Ca] J.W.S. Cassels, {\it An introduction to the geometry of numbers}, Classics in Math., Springer 1971.
\medskip
\noindent
[CMP] L. Cerlienco, M. Mignotte and F. Piras, {\it Suites r\'ecurrentes lin\'eaires: propri\'et\'es alg\'ebriques et arithm\'etiques}, L'Enseignement Math\'ematique {\bf 33} (1987), 67-108.
\medskip
\noindent
[D] H. Derksen, {\it A Skolem-Mahler-Lech theorem in positive characteristic and finite automata}, Invent. Math. {\bf 168} (2007), 175-224.  
\medskip
\noindent
[E] J.-H. Evertse, {\it On sums of S-units and linear recurrences}, Compositio Mathematica {\bf 53} (1984), 225-244.
\medskip
\noindent
[ESS] J.-H. Evertse, H.P. Schlickewei and W.M. Schmidt, {\it Linear equations in variables which lie in a multiplicative group}, Annals of Math. {\bf 155} (2002), 807-836.
\medskip
\noindent
[EZ] J.-H. Evertse and U. Zannier, {\it Linear equations with unknowns from a multiplicative group in a function field}, Acta Arith. {\bf 133} (2008), 159-170.
\medskip
\noindent
[Gh] D. Ghioca, {\it The isotrivial case in the Mordell-Lang theorem}, Trans. Amer. Math. Soc. {\bf 360} (2008), 3839-3856.
\medskip
\noindent
[GM] D. Ghioca and R. Moosa, {\it Division points on subvarieties of isotrivial semiabelian varieties}, Internat. Math. Res. Notices {\bf 19}, 2006, Article ID 65437, 1-23. 
\medskip
\noindent
[Go] D. Goss, {\it Basic structures of function field arithmetic}, Ergebnisse der Math. {\bf 35}, Springer 1996.
\medskip
\noindent
[He] D.R. Heath-Brown, {\it Artin's conjecture for primitive roots}, Quart. J. Math. Oxford Ser. (2) {\bf 37} (1986), 27-38.
\medskip
\noindent
[HP] W.V.D. Hodge and D. Pedoe, {\it Methods of algebraic geometry I}, Cambridge 1968.
\medskip
\noindent
[Hr] E. Hrushovsky, {\it The Mordell-Lang conjecture for function fields}, J. Amer. Math. Soc. {\bf 9} (1996), 667-690.
\medskip
\noindent
[HW] L.-C. Hsia and J. T.-Y. Wang, {\it The ABC theorem for higher-dimensional function fields}, Trans. Amer. Math. Soc. {\bf 356} (2003), 2871-2887.
\medskip
\noindent
[La1] S. Lang, {\it Introduction to algebraic geometry}, Addison-Wesley 1973.
\medskip
\noindent
[La2] S. Lang, {\it Fundamentals of diophantine geometry}, Springer 1983.
\medskip
\noindent
[Le] D. Leitner, {\it Linear equations in positive characteristic}, Master Thesis, University of Basel 2008.
\medskip
\noindent
[Maso] R.C. Mason, {\it Diophantine equations over function fields}, London Math. Soc. Lecture Notes {\bf 96}, Cambridge 1984.
\medskip
\noindent
[Mass] D. Masser, {\it Mixing and linear equations over groups in positive characteristic}, Israel J. Math. {\bf 142} (2004), 189-204.
\medskip
\noindent
[MS1]  R. Moosa and T. Scanlon, {\it The Mordell-Lang conjecture in positive characteristic revisited},
Model theory and applications, Quaderni di matematica {\bf 11} (eds. L. Belair, Z. Chatzidakis, P. D'Aquino, D. Marker, M. Otero, F. Point, and A. Wilkie), Dipartimento di Matematica
Seconda Universit\`a di Napoli (2002) pp. 273-296. 
\medskip
\noindent
[MS2] R. Moosa and T. Scanlon, {\it $F$-structures and integral points on semiabelian varieties over finite fields}, American Journal of Mathematics {\bf 126} (2004), 473-522.
\medskip
\noindent
[PS] A.J. van der Poorten and H. P. Schlickewei, {\it Additive relations in fields}, Journal of the Australian Mathematical Society {\bf A 51} (1991), 154-170.
\medskip
\noindent
[S] W.M. Schmidt, {\it Diophantine approximations and diophantine equations}, Lecture Notes in Math. {\bf 1467}, Springer 1991.
\medskip
\noindent
[SV] T. Struppeck and J.D. Vaaler, {\it Inequalities for heights of algebraic subspaces and the Thue-Siegel principle}, in Analytic Number Theory (Allerton Park 1989), Progress in Math. {\bf 85}, Birkh\"auser Boston 1990 (pp. 493-528).
\medskip
\noindent
[T] J. Thunder, {\it Siegel's lemma for function fields}, Michigan Math. J. {\bf 42} (1995), 147-162.
\medskip
\noindent
[V] J.F. Voloch, {\it The equation $ax+by=1$ in characteristic $p$}, J. Number Theory {\bf 73} (1998), 195-200.
\vglue 2cm
\noindent
{\bf H. Derksen:} Department of Mathematics, University of Michigan, East Hall 530 Church Street, Ann Arbor, Michigan 48104, U.S.A. ({\it hderksen@umich.edu}) 
\medskip
\noindent
{\bf D. Masser:} Mathematisches Institut, Universit\"at Basel, Rheinsprung 21, 4051 Basel, Switzerland ({\it David.Masser@unibas.ch}).

\vglue 2cm
\line {\hfil 23rd September 2010}

\end